%% file: SpaceSelecFinal4.tex
\documentclass[12pt,a4paper,final]{article}
\usepackage[a4paper,portrait]{geometry}
\usepackage{amsmath,amsfonts,latexsym,amssymb,amsthm}
\usepackage{graphicx}

\input{macros}
\THMEN

\numberwithin{equation}{section}

\begin{document}

\title{Nonparametric estimation of covariance functions by model selection}
\author{J. Bigot, R. Biscay, J-M. Loubes and L. Muniz}
\date{}
\maketitle

\begin{abstract}
We propose a model selection approach for covariance estimation of a multi-dimensional stochastic process. Under very general assumptions, observing i.i.d replications of the process at fixed observation points, we construct an estimator of the covariance function by expanding the process onto a collection of basis functions. We study the non asymptotic property of this estimate and give a tractable way of selecting the  best estimator among a possible set of candidates. The optimality of the procedure is proved via an oracle inequality which warrants that the best model is selected.
\end{abstract}


\noindent {\ \noindent \textbf{Keywords}: Model Selection, Covariance
Estimation, Stochastic process, Basis expansion, Oracle inequality. \newline
\textbf{Subject Class. MSC-2000}: 62G05, 62G20 .}

\section{Introduction}

Covariance estimation is a fundamental issue in inference for stochastic
processes with many applications, ranging from hydroscience, geostatistics,
financial series or epidemiology for instance (we refer to \cite{092462100},
\cite{MR0456314} or \cite{MR1239641} for general references for
applications). Parametric methods have been extensively studied in the
statistical literature (see \cite{MR1239641} for a review) while nonparametric procedure have received a growing attention along the last
decades. One of the main issue in this framework is to impose that the
estimator is also a covariance function, preventing the direct use of usual
nonparametric statistical methods. In this paper, we propose to use a model selection
procedure to construct a nonparametric estimator of the covariance
function of a stochastic process under general assumptions for the process. In particular we will not assume Gaussianity nor stationarity.

Consider a stochastic process $X(t)$ with values in $\mathbb{R}$, indexed by
$t\in T$, a subset of $\mathbb{R}^{d}$, $d\in \mathbb{N}$. Throughout the
paper, we assume that $X$ has finite covariance $\sigma \left( s,t\right)
=cov\left( X\left( s\right) ,X\left( t\right) \right) <+\infty $ for all $%
s,t\in T$ and, for sake of simplicity, zero mean $\mathbb{E}\left( X\left(
t\right) \right) =0$ for all $t\in T$. The observations are $X_{i}\left( t_{j}\right) $ for $i=1,...,N$, $j=1,...,n$, where
the observation points $t_{1},...,t_{n}\in T$ are fixed, and $X_{1},...,X_{N}$ are
independent copies of the process $X$. Our aim is to build a nonparametric estimator of its covariance.\\ \indent Functional approximations
of the processes $X_{1}$,...,$X_{N}$ from data $\left( X_{i}(t_{j}) \right) $ are involved in covariance function estimation. When
dealing with functional data analysis (see, e.g., \cite{Ramsey-Silverman05}), smoothing the processes $X_{1}$,...,$X_{N}$ is sometimes carried out as a
first step before computing the empirical covariance such as spline interpolation for example (see for instance in \cite{Elogne}) or projection onto a general finite basis. Let $\mathbf{x}_{i}=\left( X_{i}\left( t_{1}\right) ,...,X_{i}\left(
t_{n}\right) \right) ^{T}$ be the vector of observations at the points $%
t_{1},...,t_{n}$ with $i\in \left\{ 1,...,N\right\} .$ Let $\left\{ g_{\lambda} \right\}_{\lambda \in \mathcal{M}}$ be a collection of  (usually linearly independent but not always) functions $g_{\lambda}:T \rightarrow \mathbb{R}$ where $\mathcal{M}$ denote a generic countable set of indices. Then, let $(m) \subset \mathcal{M}$ be a subset of indices of size $m \in \mathbb{N}$ and define the $n\times m$  matrix $\mathbf{G}$ with
entries $g_{j\lambda}=g_{\lambda}\left( t_{j}\right) $, $j=1,...,n$, $\lambda \in (m)$.  $\mathbf{G}$ will be called the design matrix corresponding to the set of basis functions indexed by $(m)$. 

In such setting, usual covariance estimation is a two-step procedure: first, for each $%
i=1,...,N$, fit the regression model
\begin{equation}
\mathbf{x}_{i}=\mathbf{Ga}_{i}+\mathbf{\epsilon }_{i}  \label{RegModelXi}
\end{equation}%
(by least squares or regularized least squares), where $\mathbf{\epsilon }%
_{i}$ are random vectors in $\mathbb{R}^{n}$, to obtain estimates $\widehat{%
\mathbf{a}}_{i} = (\hat{a}_{i,\lambda})_{\lambda \in (m)} \in \mathbb{R}^{m}$ of $\mathbf{a}_{i}$ where in the case of standard least squares estimation (assuming for simplicity that $\mathbf{G}^{T}\mathbf{G}$ is invertible)
$$
\widehat{\mathbf{a}}_{i} = (\mathbf{G}^{T}\mathbf{G})^{-1} \mathbf{G}^{T} \mathbf{x}_{i}, i=1,\ldots,N.
$$
Then, estimation of  the covariance is given by computing
the following estimate
\begin{equation} \label{eq:sigma}
\widehat{\mathbf{\Sigma }}=\mathbf{G}\widehat{\mathbf{\Psi }}\mathbf{G}^{T},
\end{equation}%
where
\begin{equation}  \label{eq:psi}
\widehat{\mathbf{\Psi }}=\frac{1}{N}\sum_{i=1}^{N}\widehat{\mathbf{a}}_{i}%
\widehat{\mathbf{a}}_{i}^{T} =   (\mathbf{G}^{T}\mathbf{G})^{-1} \mathbf{G}^{T} \left( \frac{1}{N}\sum_{i=1}^{N} \mathbf{x}_{i}\mathbf{x}_{i}^T \right) \mathbf{G}  (\mathbf{G}^{T}\mathbf{G})^{-1}.
\end{equation}%
This corresponds to approximate the process $X_i$ by a truncated process $\tilde{X}_{i}$ defined as
$$
\widetilde{X}_{i}\left( t\right) =\sum\limits_{\lambda \in
(m)} \hat{a}_{i,\lambda}g_{\lambda}\left( t\right), i=1,\ldots,N,
$$
and to choose the empirical covariance of $\tilde{X}$ as an estimator of the covariance of $X$, defined by
\begin{equation*}
\widehat{\sigma }\left( s,t\right) =\frac{1}{N}\sum_{i=1}^{N}\widetilde{X_{i}%
}\left( s\right) \widetilde{X_{i}}\left( t\right) .
\end{equation*}

In this paper we propose to view the estimator (\ref{eq:sigma}) as the covariance obtained by considering a least squares estimator in the following matrix regression model
\begin{equation}
\mathbf{x}_{i}\mathbf{x}_{i}^{T}=\mathbf{G\Psi G}^{T}+\mathbf{U}_{i},\quad
i=1,...,N,  \label{MatrixRegressionModel}
\end{equation}%
where $\mathbf{\Psi }$ is a symmetric matrix and $\mathbf{U}_{i}$ are i.i.d
matrix errors.  Fitting the models (\ref{RegModelXi}) and (\ref{MatrixRegressionModel}) by least squares naturally  leads to the definition of different contrast and risk functions as the estimation is not performed in the same space ($\mathbb{R}^{m}$ for model (\ref{RegModelXi}) and $\mathbb{R}^{m \times m}$ for  model  (\ref{MatrixRegressionModel})). By choosing an appropriate loss function, least squares estimation in
model  (\ref{MatrixRegressionModel}) also leads to the natural estimate (\ref{eq:sigma}) derived from least square estimation in model (\ref{RegModelXi}). However, the
problem of model selection, i.e. choosing an appropriate data-based subset of  indices $(m) \in \mathcal{M}$, is very distinct in model (\ref{RegModelXi}) and  model (\ref%
{MatrixRegressionModel}). Indeed, model selection for (\ref{RegModelXi})
depends on the variability of the vectors $\mathbf{x}_{i}$'s while for (\ref%
{MatrixRegressionModel}) it depends on the variability of the matrices $\mathbf{x}_{i}%
\mathbf{x}_{i}^{T}$'s. One of the main contributions of this paper is to show that considering model (\ref{MatrixRegressionModel}) enables to handle a large variety of cases and to build an optimal model selection estimator of the covariance without too strong assumptions on the model. Moreover it will be shown that considering model  (\ref{MatrixRegressionModel}) leads to the estimator  $\widehat{\mathbf{\Psi }}$ (\ref{eq:psi}) which is guaranteed to be in the class of definite non negative matrices and thus to  a proper covariance matrix $\widehat{\mathbf{\Sigma }} = \mathbf{G}\widehat{\mathbf{\Psi }}\mathbf{G}^{T} $.

A similar method has been developed for smooth
interpolation of covariance functions in \cite{Biscay-etal05}. However, this paper is restricted to basis functions that are determined by
reproducing kernels in suitable Hilbert spaces. Furthermore, a matrix metric
different from (though related to) the Frobenius matrix norm is adopted as a
fitting criterion.  Similar ideas are tackled in \cite{Matsuo-etal06}. These
authors deal with the estimation of $\mathbf{\Sigma }$ within the
covariance class $\mathbf{\Gamma }=\mathbf{G\Psi G}^{T}$ induced by an
orthogonal wavelet expansion. However, their fitting criterion is not general since they choose the Gaussian
likelihood as a contrast function, and thus their method requires specific
distributional assumptions. We also point out that computation of the Gaussian likelihood requires
inversion of $\mathbf{G\Psi G}^{T}$, which is not directly feasible if $%
rank\left( \mathbf{G}\right) <n$ or some diagonal entities of the definite non negative (d.n.n)
matrix $\mathbf{\Psi }$ are zero. \\ \indent Hence, to our knowledge, no previous work has proposed to use the matrix regression
model (\ref{MatrixRegressionModel}) under general moments assumptions of the process $X$ using a general basis expansion for nonparametric covariance function estimation.\vskip .1in The paper then falls into the following parts. The description of the
statistical framework of the matrix regression is given in Section~\ref{s1}. Section 2 is devoted to the main statistical results. Namely we study the behavior of the estimator for a fixed model in Section~2.1 while Section~2.2 deals with the model selection procedure and provide the oracle inequality. Section~3 states a concentration inequality that is used in all the paper, while the proofs are postponed to
a technical Appendix .

\section{Nonparametric Model selection for Covariance estimation} \label{s1}

Recall that $X=\left( X\left( t\right) \right) _{t\in T}$ is an $\mathbb{R}$%
-valued stochastic process, where $T$ denotes some subset of $\mathbb{R}^{d}$%
, $d\in \mathbb{N}$. Assume that $X$ has finite moments up to order $4$, and
zero mean, i.e $\mathbb{E}\left( X\left( t\right) \right) =0$ for all $t\in
T $. The covariance function of $X$ is denoted by $\sigma \left( s,t\right)
=cov\left( X\left( s\right) ,X\left( t\right) \right) $ for $s,t\in T$ and recall that  $%
X_{1},...,X_{N}$ are independent copies of the process $X$.\newline
\indent In this work, we observe at different observation points $%
t_{1},...,t_{n}\in T$ these independent copies of the process, denoted by $%
X_{i}\left( t_{j}\right) $, with $i=1,...,N$, $j=1,...,n$. Recall that $\mathbf{x}%
_{i}\mathbf{=}\left( X_{i}\left( t_{1}\right) ,...,X_{i}\left( t_{n}\right)
\right) ^{T}$ is the vector of observations at the points $t_{1},...,t_{n}$
for each $i=1,...,N$. The matrix $\mathbf{\Sigma =}\mathbb{E}\left( \mathbf{x%
}_{i}\mathbf{x}_{i}^{T}\right) =\left( \sigma \left( t_{j},t_{k}\right)
\right) _{1\leq j\leq n,1\leq k\leq n}$ is the covariance matrix of $X$ at the observations points. Let
$\overline{\mathbf{x}}$ and $\mathbf{S}$ denote the sample mean and the
sample covariance (non corrected by the mean) of the data $\mathbf{x}_{1},...,%
\mathbf{x}_{N}$, i.e.
\begin{equation*}
\overline{\mathbf{x}}\mathbf{=}\frac{1}{N}\sum_{i=1}^{N}\mathbf{x}%
_{i},\qquad \mathbf{S=}\frac{1}{N}\sum_{i=1}^{N}\mathbf{x}_{i}\mathbf{x}%
_{i}^{T}.
\end{equation*}
Our aim is to build a model selection estimator of the covariance of the
process observed with $N$ replications but without additional assumptions
such as stationarity nor Gaussianity. The asymptotics will be taken with
respect to $N$, the number of copies of the process.

\subsection{Notations and preliminary definitions}

First, define specific matricial notations. We refer to \cite{Lutkepohl96}
or \cite{MR2162145} for definitions and properties of matrix operations and
special matrices.  As usual, vectors in $\mathbb{R}^{k}$ are regarded as column vectors
for all $k\in \mathbb{N}$. To be able to write general methods for all our
models, we will treat matricial data as a natural extension of the vectorial
data, with of course, different correlation structure. For this, we
introduce a natural linear transformation, which converts any matrix into a
column vector. The vectorization of a $k\times n$ matrix $\mathbf{A}%
=(a_{ij})_{1\leq i\leq k,1\leq j\leq n}$ is the $kn\times 1$ column vector
denoted by $vec\left( \mathbf{A} \right) $, obtain by stacking the columns
of the matrix $\mathbf{A}$ on top of one another. That is $%
vec(A)=[a_{11},...,a_{k1},a_{12},...,a_{k2},...,a_{1n},...,a_{kn}]^{T}$.%
\newline
\indent For a symmetric $k\times k$ matrix $\mathbf{A}$, the vector $%
vec\left( \mathbf{A}\right) $ contains more information than necessary,
since the matrix is completely determined by the lower triangular portion,
that is, the $k(k+1)/2$ entries on and below the main diagonal. Hence, we
introduce the symmetrized vectorization, which corresponds to a
half-vectorization, denoted by $vech(\mathbf{A})$. More precisely, for any
matrix $\mathbf{A}=(a_{ij})_{1\leq i\leq k,1\leq j\leq k}$, define $vech(%
\mathbf{A})$ as the $k(k+1)/2\times 1$ column vector obtained by vectorizing
only the lower triangular part of $\mathbf{A}$. That is $vech(\mathbf{A}%
)=[a_{11},...,a_{k1},a_{22},...,a_{n2},...,a_{\left( k-1\right) \left(
k-1\right) },a_{\left( k-1\right) k},a_{kk}]^{T}$. There exist unique linear
transformation which transforms the half-vectorization of a matrix to its
vectorization and vice-versa called, respectively, the duplication matrix
and the elimination matrix. For any $k\in \mathbb{N}$, the $k^{2}\times
k\left( k+1\right) /2$ duplication matrix is denoted by $\mathbf{D}_{k}$, $%
\mathbf{1}_{k}=\left( 1,...,1\right) ^{T}\in \mathbb{R}^{k}$ and $\mathbf{I}%
_{k}$ is the identity matrix in $\mathbb{R}^{k\times k}$. 

For any matrix $\mathbf{A}$, $\mathbf{A}^{T}$ is the transpose of $\mathbf{A}
$\textbf{, }$tr\left( \mathbf{A}\right) $ is the trace of $\mathbf{A}$, $%
\left\Vert \mathbf{A}\right\Vert $ is the Frobenius matrix norm defined as $%
\left\Vert \mathbf{A}\right\Vert ^{2}=tr\left( \mathbf{AA}^{T}\right) $, $%
\lambda _{\max }\left( \mathbf{A}\right) $ is the maximum eigenvalue of $%
\mathbf{A}$\textbf{, }$\rho \left( \mathbf{A}\right) $ is the spectral norm
of $\mathbf{A}$, that is $\rho \left( \mathbf{A}\right) =\lambda _{\max
}\left( \mathbf{A}\right) $ for $\mathbf{A}$ a d.n.n
matrix. If $\mathbf{A=}(a_{ij})_{1\leq i\leq k,1\leq j\leq n}$ is a $k\times n$
matrix and $\mathbf{B=}(b_{ij})_{1\leq i\leq p,1\leq j\leq q}$ is a $p\times
q$ matrix, then the Kronecker product of the two matrices, denoted by $%
\mathbf{A}\otimes \mathbf{B}$, is the $kp\times nq$ block matrix%
\begin{equation*}
\mathbf{A}\otimes \mathbf{B=}%
\begin{bmatrix}
a_{11}\mathbf{B} & . & . & . & a_{1n}\mathbf{B} \\
. & . &  &  & . \\
. &  & . &  & . \\
. &  &  & . & . \\
a_{k1}\mathbf{B} & . & . & . & a_{kn}\mathbf{B}%
\end{bmatrix}%
.
\end{equation*}

For any random matrix $\mathbf{Z=}\left( Z_{ij}\right) _{1\leq i\leq k,1\leq
j\leq n}$, its expectation is denoted by $\mathbb{E}\left( \mathbf{Z}\right)
=\left( \mathbb{E}\left( Z_{ij}\right) \right) _{1\leq i\leq k,1\leq j\leq n}
$. For any random vector $\mathbf{z=}\left( Z_{i}\right) _{1\leq i\leq k}$,
let $V\left( \mathbf{z}\right) =\left( cov\left( Z_{i},Z_{j}\right) \right)
_{1\leq i,j\leq k}$ be its covariance matrix. With this notation, $V\left(
\mathbf{x}_{1}\right) =V\left( \mathbf{x}_{i}\right) =\left( \sigma \left(
t_{j},t_{k}\right) \right) _{1\leq j\leq n,1\leq k\leq n}$ is the covariance
matrix of $X$.

Let $(m) \in \mathcal{M}$, and recall that to the finite set $\mathcal{G}_m=\left\{ g_{\lambda} \right\}_{\lambda \in (m)} $ of  functions $g_{\lambda}:T\rightarrow%
\mathbb{R}$ we associate the $n\times m$ matrix $\mathbf{G}$ with
entries $g_{j\lambda}=g_{\lambda}\left( t_{j}\right) $, $j=1,...,n$, $\lambda \in (m)$. Furthermore, for each $t\in T$, we write $%
\mathbf{G}_{t}=\left( g_{\lambda}\left( t\right), \lambda \in (m)  \right)
^{T}$. For $k\in \mathbb{N}$, $\mathcal{S}_{k}$ denotes the linear subspace of $%
\mathbb{R}^{k\times k}$ composed of symmetric matrices. For $\mathbf{G\in }%
\mathbb{R}^{n\times m}$, $\mathcal{S}\left( \mathbf{G}\right) $ is the
linear subspace of $\mathbb{R}^{n\times n}$ defined by%
\begin{equation*}
\mathcal{S}\left( \mathbf{G}\right) =\left\{ \mathbf{G\Psi G}^{T}:\mathbf{%
\Psi \in }\mathcal{S}_{m}\right\} \text{.}
\end{equation*}
Let $\mathcal{S}_{N}\left( \mathbf{G}\right) $ be the linear subspace of $%
\mathbb{R}^{nN\times n}$ defined by

\begin{equation*}
\mathcal{S}_{N}\left( \mathbf{G}\right) =\left\{ \mathbf{1}_{N}\otimes
\mathbf{G\Psi G}^{T}:\mathbf{\Psi \in }\mathcal{S}_{m}\right\} =\left\{
\mathbf{1}_{N}\otimes \mathbf{\Gamma }:\mathbf{\Gamma \in }\mathcal{S}\left(
\mathbf{G}\right) \right\}
\end{equation*}%
and let $\mathcal{V}_{N}\left( \mathbf{G}\right) $ be the linear subspace of
$\mathbb{R}^{n^{2}N}$ defined by%
\begin{equation*}
\mathcal{V}_{N}\left( \mathbf{G}\right) =\left\{ \mathbf{1}_{N}\otimes
vec\left( \mathbf{G\Psi G}^{T}\right) :\mathbf{\Psi \in }\mathcal{S}%
_{m}\right\} =\left\{ \mathbf{1}_{N}\otimes vec\left( \mathbf{\Gamma }%
\right) :\mathbf{\Gamma \in }\mathcal{S}\left( \mathbf{G}\right) \right\}
\text{.}
\end{equation*}%
All these spaces are regarded as Euclidean spaces with the scalar product
associated to the Frobenius matrix norm.

\subsection{Model}

The approach that we will develop to estimate the covariance function $%
\sigma $ is based on the following two main ingredients: first,
we consider a functional expansion $\tilde{X}$ to approximate the underlying process $X$
and take the covariance of $\tilde{X}$  as an approximation of the true
covariance $\Sigma$.\newline
\indent For this, let $(m) \in \mathcal{M}$ and consider an approximation to the process $X$ of the following form:
\begin{equation}
\widetilde{X}\left( t\right) =\sum\limits_{\lambda \in
(m)}a_{\lambda}g_{\lambda}\left( t\right) ,  \label{Expan_X}
\end{equation}
where $a_{\lambda}$ are suitable random coefficients. For instance if $X$ takes its values in $L^{2}(T)$ (the space  of square integrable real-valued functions on $T$) and if $(g_{\lambda})_{\lambda \in \mathcal{M}}$ are orthonormal  functions in $L^{2}(T)$, then one can take
$$
a_{\lambda} = \int_{T} X(t) g_{\lambda}(t) dt.
$$
Several basis can thus be
considered, such as a polynomial basis on $\mathbb{R}^{d}$, Fourier
expansion on a rectangle $T\subset\mathbb{R}^{d}$ (i.e. $g_{\lambda}\left( t\right)
=e^{i2\pi\left\langle \omega_{\lambda},t\right\rangle }$, using a regular grid of
discrete set of frequencies $\left\{ \omega_{\lambda} \in\mathbb{R}^{d}, \lambda \in (m) \right\}$ that
do not depend on $t_{1},...,t_{n}$). One can also use, as in \cite{Elogne}%
, tensorial product of B-splines on a rectangle $T\subset\mathbb{R}^{d}$,
with a regular grid of nodes in $\mathbb{R}^{d}$ not
depending on $t_{1},...,t_{n}$ or a standard wavelet basis on $\mathbb{R}^{d}$, depending on a regular grid of locations in $ \mathbb{R}%
^{d} $ and discrete scales in $\mathbb{R}_{+}$. Another
class of natural expansion is provided by Karhunen-Loeve expansion of the
process $X$ (see \cite{adler} for more references). 

Therefore, it is natural to consider the covariance function $\rho $ of $%
\widetilde{X}$ as an approximation of $\sigma $. Since the covariance $\rho $
can be written as
\begin{equation}
\rho \left( s,t\right) =\mathbf{G}_{s}^{T} \mathbf{\overline{\Psi} G}_{t},
\label{Expan_s}
\end{equation}%
where, after reindexing the functions if necessary, $\mathbf{G}_{t}=\left( g_{\lambda}\left( t\right), \lambda \in (m)  \right)
^{T}$ and
$$
\mathbf{%
\overline{\Psi} =}\left( \mathbb{E}\left( a_{\lambda}a_{\mu}\right) \right), \mbox{ with } (\lambda,\mu) \in (m) \times (m).
$$
Hence we are led to look for an estimate $\widehat{\sigma }$ of $%
\sigma $ in the class of functions of the form \eqref{Expan_s}, with $%
\mathbf{\Psi }\in \mathbb{R}^{m\times m}$ some symmetric
matrix. Note that the choice of the function expansion in \eqref{Expan_X},
in particular the choice the subset of indices $(m)$, will be crucial in
the approximation properties of the covariance function $\rho $. This
estimation procedure has several advantages: it will be shown that an appropriate choice of loss function leads to the construction of symmetric d.n.n matrix $\widehat{\mathbf{\Psi }}$ (see Proposition \ref{pro:preliminar}) and thus the resulting estimate
\begin{equation*}
\widehat{\sigma }\left( s,t\right) =\mathbf{G}_{s}^{T}\widehat{\mathbf{\Psi }}\mathbf{G}_{t},
\end{equation*}%
is a covariance function, so the resulting estimator can be plugged
in other procedures which requires working with a covariance function. We
also point out that the large amount of existing approaches for function
approximation of the type \eqref{Expan_X} (such as those based on Fourier,
wavelets, kernel, splines or radial functions) provides great flexibility to
the model \eqref{Expan_s}.

Secondly, we use the Frobenius
matrix norm to quantify the risk of the covariance matrix estimators. Recall that $\mathbf{\Sigma =}\left( \sigma \left(
t_{j},t_{k}\right) \right) _{1\leq j,k\leq n}$ is the true covariance while
$\mathbf{\Gamma =}\left( \rho \left( t_{j},t_{k}\right) \right) _{(j,k)}$ will denote
be the covariance matrix of the approximated process $\tilde{X}$ at the observation points. Hence
\begin{equation}
\mathbf{\Gamma }=\mathbf{G \overline{\Psi} G}^{T}\text{.}  \label{InducedModelSigma}
\end{equation}%
Comparing the covariance function $\rho $ with the true one $\sigma $ over
the design points $t_{j}$, implies quantifying the deviation of $\mathbf{%
\Gamma }$ from $\mathbf{\Sigma }$. For this consider the following loss
function
\begin{equation*}
L\left( \mathbf{\Psi }\right) =\mathbb{E}\left\Vert \mathbf{x}\mathbf{x}^{T}-\mathbf{G\Psi G}^{T}\right\Vert ^{2},
\end{equation*}%
where $\mathbf{x}\mathbf{=}\left( X\left( t_{1}\right)
,...,X\left( t_{n}\right) \right) ^{T}$ and $\left\Vert {.}\right\Vert $
is the Frobenius matrix norm. Note that
\begin{equation*}
L\left( \mathbf{\Psi }\right) =\left\Vert \mathbf{\Sigma }-\mathbf{G\Psi G}%
^{T}\right\Vert ^{2}+C,
\end{equation*}%
where the constant $C$ does not depend on $\mathbf{\Psi }$. The Frobenius
matrix norm provides a meaningful metric for comparing covariance matrices,
widely used in multivariate analysis, in particular in the theory on
principal components analysis. See also \cite{Biscay-etal97}, \cite%
{Schafer-Strimmer05} and references therein for other applications of this
loss function.

 To the loss $L$ corresponds the following empirical contrast
function $L_{N}$, which will be the fitting criterion we will try to minimize
\begin{equation*}
L_{N}\left( \mathbf{\Psi }\right) =\frac{1}{N}\sum_{i=1}^{N}\left\Vert
\mathbf{x}_{i}\mathbf{x}_{i}^{T}-\mathbf{G\Psi G}^{T}\right\Vert ^{2}.
\end{equation*}%
We point out that this loss is exactly the sum of the squares of the
residuals corresponding to the matrix linear regression model
\begin{equation}
\mathbf{x}_{i}\mathbf{x}_{i}^{T}=\mathbf{G\Psi G}^{T}+\mathbf{U}_{i},\quad
i=1,...,N,  \label{MatrixRegModel}
\end{equation}%
with i.i.d. matrix errors $\mathbf{U}_{i}$ such that $\mathbb{E}\left( \mathbf{U}%
_{i}\right) =\mathbf{0}$. This remark provides a natural framework to study
the covariance estimation problem as a matricial regression model. Note also
that the set of matrices $\mathbf{G\Psi G}^{T}$ is a linear subspace of $%
\mathbb{R}^{n\times n}$ when $\mathbf{\Psi }$ ranges over the space of
symmetric matrices $\mathcal{S}_{m}$.

To summarize our approach, we finally propose following two-step estimation procedure: in a first step, for a given design matrix $\mathbf{G}$, define
\begin{equation*}
\widehat{\mathbf{\Psi }}=\mathrm{arg}\min_{\mathbf{\Psi }\in \mathcal{S}%
_{m}}L_{N}(\mathbf{\Psi }),
\end{equation*}%
and take $\widehat{\mathbf{\Sigma }}=\mathbf{G}\widehat{\mathbf{\Psi }}%
\mathbf{G}^{T}$ as an estimator of $\mathbf{\Sigma }$. Note that $%
\widehat{\mathbf{\Psi }}$ will be shown to be a d.n.n matrix  (see Proposition \ref{pro:preliminar}) and thus $\widehat{%
\mathbf{\Sigma }}$ is also a d.n.n matrix. Since the minimization of $%
L_{N}\left( \mathbf{\Psi }\right) $ with respect to $\mathbf{\Psi }$ is done over the linear space 
of symmetric matrices $\mathcal{S}_{m}$, it can be transformed to a classical least squares linear
problem, and the computation of  $\widehat{\mathbf{\Psi }}$ is therefore quite simple. For a given design matrix $\mathbf{%
G}$, we will construct an estimator for\textbf{\ }$\Gamma =\mathbf{G  \overline{\Psi} G}%
^{T}$ which will be close to $\mathbf{\Sigma }=V\left(\mathbf{x}_{1}\right)$ as soon as $\tilde{X}$
is a sharp estimation of $X$. So, the role of $\mathbf{G}$ and thus the choice of the subset of indices $(m)$ is crucial since
it determines the behavior of the estimator.

 Hence, in second step, we aim at selecting the best design matrix $\mathbf{G} = \mathbf{G}_{m}$ among a collection
of candidates $\left\{ \mathbf{G}_{m}, (m)\in \mathcal{M} \right\}$. For this, methods and
results from the theory of model selection in linear regression can be
applied to the present context. In particular the results in \cite{Baraud00}%
, \cite{Comte} or \cite{loublud1} will be useful in dealing with model
selection for the framework \eqref{MatrixRegModel}. Note that only
assumptions about moments, not specific distributions of the data, are
involved in the estimation procedure.

\begin{rem}
We consider here a least-squares estimates of the covariance. Note that suitable regularization terms or constraints could also be
incorporated into the minimization of $L_{N}\left( \mathbf{\Psi}\right) $ to
impose desired properties for the resulting estimator, such as smoothness or
sparsity conditions as in \cite{1137.62338}. \end{rem}

\section{Oracle inequality for Covariance Estimation}\label{s:covest} 
The first part of this section describes the properties of
the least squares estimator $\widehat{\mathbf{\Sigma }}=\mathbf{G}\widehat{%
\mathbf{\Psi }}\mathbf{G}^{T}$ while the second part builds a selection
procedure to pick automatically the best estimate among a collection of
candidates.

\subsection{Least Squares Covariance Estimation}

Given some $n\times m$ fixed design matrix $\mathbf{G}$ associated to a
finite family of $m$ basis functions, the least squares covariance estimator
of $\mathbf{\Sigma}$ is defined by
\begin{equation}
\widehat{\mathbf{\Sigma}}=\mathbf{G}\widehat{\mathbf{\Psi}}\mathbf{G}%
^{T}=\arg\min\left\{ \frac{1}{N}\sum_{i=1}^{N}\left\Vert \mathbf{x}_{i}%
\mathbf{x}_{i}^{T}-\mathbf{\Gamma}\right\Vert ^{2}:\mathbf{\Gamma=G\Psi G}%
^{T},\mathbf{\Psi\in}\mathcal{S}_{m}\right\} .  \label{DefLSEst}
\end{equation}
The corresponding estimator of the covariance function $\sigma$ is%
\begin{equation}
\widehat{\sigma}\left( s,t\right) =\mathbf{G}_{s}^{T}\widehat{\mathbf{\Psi}}%
\mathbf{G}_{t}.  \label{sigmaTecho}
\end{equation}

\begin{prop} \label{pro:preliminar}
\label{LemaLS-solution}Let $\mathbf{Y}_{1},...,\mathbf{Y}_{N}\in
\mathbb{R}^{n\times n}$ and $\mathbf{G\in}\mathbb{R}^{n\times m}$ be arbitrary
matrices Then, the infimum%
\[
\inf\left\{  \frac{1}{N}\sum_{i=1}^{N}\left\Vert \mathbf{Y}_{i}-\mathbf{G\Psi
G}^{T}\right\Vert ^{2}:\mathbf{\Psi\in}\mathcal{S}_{m}\right\}
\]
is achieved at
\begin{equation}
\widehat{\mathbf{\Psi}}=\left(  \mathbf{G}^{T}\mathbf{G}\right)
^{-}\mathbf{G}^{T}\left(  \frac{\overline{\mathbf{Y}}+\overline{\mathbf{Y}%
}^{T}}{2}\right)  \mathbf{G}\left(  \mathbf{G}^{T}\mathbf{G}\right)
^{-},\label{GeneralPhiTecho}%
\end{equation}
where $\left(  \mathbf{G}^{T}\mathbf{G}\right)  ^{-}$ is any generalized
inverse of $\mathbf{G}^{T}\mathbf{G}$ (see \cite{engl96} for a general definition), and
\[
\overline{\mathbf{Y}}\mathbf{=}\frac{1}{N}\sum_{i=1}^{N}\mathbf{Y}%
_{i}\mathbf{.}%
\]
Furthermore, $\mathbf{G\widehat{\mathbf{\Psi}}G}^{T}$ is the same for all the
generalized inverses $\left(  \mathbf{G}^{T}\mathbf{G}\right)  ^{-}$ of
$\mathbf{G}^{T}\mathbf{G}$. In particular, if $\mathbf{Y}_{1},...,\mathbf{Y}_{N}\in\mathcal{S}_{n}$
(i.e., if they are symmetric matrices) then any minimizer has the form
\[
\widehat{\mathbf{\Psi}}=\left(  \mathbf{G}^{T}\mathbf{G}\right)
^{-}\mathbf{G}^{T}\overline{\mathbf{Y}}\mathbf{G}\left(  \mathbf{G}%
^{T}\mathbf{G}\right)  ^{-}.
\]
If $\mathbf{Y}_{1},...,\mathbf{Y}_{N}$ are d.n.n. then these matrices
$\widehat{\mathbf{\Psi}}$ are d.n.n.
\end{prop}

If we assume that $(\mathbf{G}^T \mathbf{G})^{-1}$ exists, then Proposition \ref{pro:preliminar} shows that we retrieve the expression (\ref{eq:psi}) for $\widehat{\mathbf{\Psi}}$ that has been derived from least square estimation in model (\ref{RegModelXi}).


\begin{thm}
\label{LS-solution} Let $\mathbf{S=}\frac{1}{N}\sum_{i=1}^{N}\mathbf{x}_{i}\mathbf{x}%
_{i}^{T}$. Then, the least squares covariance estimate defined by
(\ref{DefLSEst}) is given by the d.n.n. matrix
\[
\widehat{\mathbf{\Sigma}}=\mathbf{G}\widehat{\mathbf{\Psi}}\mathbf{G}^{T}%
=\Pi\mathbf{S}\Pi,
\]
where%
\begin{align}
\widehat{\mathbf{\Psi}}  & =\left(  \mathbf{G}^{T}\mathbf{G}\right)
^{-}\mathbf{G}^{T}\mathbf{SG}\left(  \mathbf{G}^{T}\mathbf{G}\right)
^{-},\label{PhiTecho}\\
\Pi & =\mathbf{G}\left(  \mathbf{G}^{T}\mathbf{G}\right)  ^{-}\mathbf{G}%
^{T}.\nonumber
\end{align}
Moreover $\widehat{\mathbf{\Sigma}}$ has the following interpretations in terms of
orthogonal projections: 

$i)$ $\widehat{\mathbf{\Sigma}}$ is the projection of $\mathbf{S}\in\mathbb{R}%
^{n\times n}$ on $\mathcal{S}\left(  \mathbf{G}\right)  $.

$ii)$ $\mathbf{1}_{N}\otimes\widehat{\mathbf{\Sigma}}$ is the projection of
$\mathbf{Y}=\left(  \mathbf{x}_{1}\mathbf{x}_{1}^{T},...,\mathbf{x}%
_{N}\mathbf{x}_{N}^{T}\right)  ^{T}\in\mathbb{R}^{nN\times n}$ on
$\mathcal{S}_{N}\left(  \mathbf{G}\right)  .$

$iii)$ $\mathbf{1}_{N}\otimes vec\left(  \widehat{\mathbf{\Sigma}}\right)  $ is the
projection of $\mathbf{y}=\left(  vec^{T}\left(  \mathbf{x}_{1}\mathbf{x}%
_{1}^{T}\right)  ,...,vec^{T}\left(  \mathbf{x}_{N}\mathbf{x}_{N}^{T}\right)
\right)  ^{T}\in\mathbb{R}^{n^{2}N}$ on $\mathcal{V}_{N}\left(  \mathbf{G}%
\right)  .$

\end{thm}

The proof of this theorem is a direct application of Proposition~%
\ref{LemaLS-solution}. Hence for a given design matrix $\mathbf{G}$%
, the least squares estimator $\mathbf{\hat{\Sigma}}=\mathbf{\hat{\Sigma}}(%
\mathbf{G})$ is well defined and has the structure of a covariance matrix.
It remains to study how to pick automatically the estimate when dealing with
a collection of design matrices coming from several approximation choices
for the random process $X$.

\subsection{Main Result}

Consider a collection of indices $(m)\in \mathcal{M}$ with size $m$.
Let also $\left\{ \mathbf{G}_{m}:(m)\in \mathcal{M}\right\} $ be a
finite family of design matrices $\mathbf{G}_{m}\in \mathbb{R}^{n\times m}$,
and let $\widehat{\mathbf{\Sigma }}_{m} = \mathbf{\hat{\Sigma}}(%
\mathbf{G}_{m})$, $(m)\in \mathcal{M}$, be the
corresponding least squares covariance estimators. The problem of interest
is to select the best of these estimators in the sense of the minimal
quadratic risk $\mathbb{E}\left\Vert \mathbf{\Sigma }-\widehat{\mathbf{%
\Sigma }}_{m}\right\Vert ^{2}.$

The main theorem of this section provides a non-asymptotic bound for the
risk of a penalized strategy for this problem. For all $(m)\in \mathcal{M}$, write%
\begin{align}
\mathbf{\Pi }_{m}& =\mathbf{G}_{m}\left( \mathbf{G}_{m}^{T}\mathbf{G}%
_{m}\right) ^{-}\mathbf{G}_{m}^{T},  \label{Pi-m} \\
D_{m}& =Tr\left( \mathbf{\Pi }_{m}\right) ,  \notag
\end{align}%
We assume that $D_{m}\geq 1$ for all $(m) \in \mathcal{M}.$ The estimation error for a given model $(m) \in \mathcal{M}$ is given by
\begin{equation}
\mathbb{E}\left( \left\Vert \mathbf{\Sigma }-\widehat{\mathbf{\Sigma }}%
_{m}\right\Vert ^{2}\right) =\left\Vert \mathbf{\Sigma -\Pi }_{m}\mathbf{%
\Sigma \Pi }_{m}\right\Vert ^{2}+\frac{\delta _{m}^{2}D_{m}}{N},
\label{simple}
\end{equation}%
where%
\begin{align*}
\delta _{m}^{2}& =\frac{\mathrm{Tr}\left( \left( \mathbf{\Pi }_{m}\otimes
\mathbf{\Pi }_{m}\right) \mathbf{\Phi }\right) }{D_{m}}, \\
\mathbf{\Phi }& \mathbf{=}V\left( vec\left( \mathbf{x}_{1}\mathbf{x}%
_{1}^{T}\right) \right) .
\end{align*}
Given $\theta >0$, define the penalized covariance estimator $\widetilde{%
\mathbf{\Sigma }}=\widehat{\mathbf{\Sigma }}_{\widehat{m}}$ by
\begin{equation*}
\widehat{m}=\arg \underset{(m)\in \mathcal{M}}{\min }\left\{ \frac{1}{N%
}\sum_{i=1}^{N}\left\Vert \mathbf{x}_{i}\mathbf{x}_{i}^{T}-\widehat{\mathbf{%
\Sigma }}_{m}\right\Vert ^{2}+pen\left( m\right) \right\} ,
\end{equation*}%
where
\begin{equation}
pen\left( m\right) =\left( 1+\theta \right) \frac{\delta _{m}^{2}D_{m}}{N}.
\label{Penal}
\end{equation}%
\begin{thm}
\label{TeorConcentIneqLSEst}
Let $q>0$ be given such that there exists $p>2\left(  1+q\right)  $ satisfying
$\mathbb{E}\left\Vert \mathbf{x}_{1}\mathbf{x}_{1}^{T}\right\Vert ^{p}<\infty
$. Then, for some constants $K\left(  \theta\right)  >1$ and $C^{\prime
}\left(  \theta,p,q\right)  >0$ we have that
\[
\left(  \mathbb{E}  \left\Vert \mathbf{\Sigma}-\widetilde{\mathbf{\Sigma
}}\right\Vert ^{2q}\right)  ^{1/q}\leq2^{\left(  q^{-1}-1\right)
_{+}}\left[  K\left(  \theta\right)  \inf_{(m)\in\mathcal{M}}\left(
\left\Vert \mathbf{\Sigma-\Pi}_{m}\mathbf{\Sigma\Pi}_{m}\right\Vert ^{2}%
+\frac{\delta_{m}^{2}D_{m}}{N}\right)  +\frac{\Delta_{p}}{N}\delta_{\sup}%
^{2}\right]  ,
\]
where
\[
\Delta_{p}^{q}=C^{\prime}\left(  \theta,p,q\right)  \mathbb{E}\left\Vert
\mathbf{x}_{1}\mathbf{x}_{1}^{T}\right\Vert ^{p}\left(  \sum\limits_{(m)\in
\mathcal{M}}\delta_{m}^{-p}D_{m}^{-\left(  p/2-1-q\right)  }\right)
\] and
$$ \delta_{\sup}^{2}   =\max\left\{  \delta_{m}^{2}: (m) \in\mathcal{M}%
\right\}  .$$
In particular, for $q=1$ we have
\begin{equation}
\mathbb{E}\left(  \left\Vert \mathbf{\Sigma}-\widetilde{\mathbf{\Sigma}%
}\right\Vert ^{2}\right)  \leq K\left(  \theta\right)  \inf_{(m) \in
\mathcal{M}} \mathbb{E}\left( \left\Vert \mathbf{\Sigma }-\widehat{\mathbf{\Sigma }}%
_{m}\right\Vert ^{2}\right)  +\frac{\Delta_{p}%
}{N}\delta_{\sup}^{2}.\label{QuadConIneq}%
\end{equation}

\end{thm}For the proof of this result, we first restate this theorem in a a
vectorized form which turns to be a $d$-variate extensions of results in
\cite{Baraud00} (which are covered when $d=1$) and are stated in Section~\ref%
{s:auxBaraud}. Their proof rely on model selection techniques and a
concentration tool stated in Section~\ref{s:concentration}.

\begin{rem}The penalty depends on the quantity $\delta_m$. Note that
\begin{align}
D_{m}\delta_{m}^{2}  & =\gamma_{m}^{2}=\gamma^{2}\left(  m,n\right)
={\rm Tr}\left(  \left(  \mathbf{\Pi}_{m}\otimes\mathbf{\Pi}_{m}\right)
\mathbf{\Phi}\right) \label{Gamma-m}\\
& =\mathbb{E}\left\Vert \widehat{\mathbf{\Sigma}}_{m}-\mathbf{\Pi}_{m}\mathbf{\Sigma
\Pi}_{m}\right\Vert ^{2}N={\rm Tr}\left(  V\left(  vec\left(  \widehat
{\mathbf{\Sigma}}_{m}\right)  \right)  \right)  N.\nonumber
\end{align}
So, we get that $\delta_{m}^{2}\leq\lambda_{\max}\left(  \mathbf{\Phi
}\right)  $ for all $(m)$. Hence Theorem~\ref{TeorConcentIneqLSEst} remains true
if $\delta_{m}^{2}$ is replaced by $\lambda^{2}=\lambda_{\max}\left(
\mathbf{\Phi}\right)  $ in all the statements.
\end{rem}

\begin{rem} The penalty relies thus on
$\mathbf{\Phi}  \mathbf{=}V\left( vec\left( \mathbf{x}_{1}\mathbf{x}%
_{1}^{T}\right) \right).$ This quantity reflects the correlation structure of the data. We point out that for practical purpose, this quantity can be estimated using the empirical version of $\Phi$ since the $\mathbf{x}_i,\: i=1,\dots,N$ are i.i.d observed random variables. In the original paper by Baraud \cite{MR1918295}, an estimator of the variance is proposed to overcome this issue. However, the consistency proof relies on a concentration inequality which turns to be a $\chi^2$ like inequality. Extending this inequality to our case would mean to be able to construct concentration bounds for matrices $\mathbf{x}\mathbf{x}^T$, implying Wishart distributions. If some results exist in this framework \cite{MR2485015}, adapting this kind of construction to our case falls beyond the scope of this paper.
\end{rem}We have obtained in Theorem~\ref{TeorConcentIneqLSEst} an oracle
inequality since, using \eqref{simple} and \eqref{QuadConIneq}, one
immediately sees that $\widetilde{\mathbf{\Sigma }}$ has the same quadratic
risk as the \textquotedblleft oracle\textquotedblright\ estimator except for
an additive term of order $O\left( \frac{1}{N}\right) $ and a constant
factor. Hence, the selection procedure is optimal in the sense that it
behaves as if the true model were at hand. To describe the result in terms
of rate of convergence, we have to pay a special attention to the bias terms
$\left\Vert \mathbf{\Sigma -\Pi }_{m}\mathbf{\Sigma \Pi }_{m}\right\Vert
^{2} $. In a very general framework, it is difficult to evaluate such
approximation terms. If the process has bounded second moments, i.e for all $%
i=1,\dots ,n$, we have $\mathbb{E}\left( X^{2}\left( t_{i}\right) \right)
\leq C$, then we can write
\begin{eqnarray*}
\left\Vert \mathbf{\Sigma -\Pi }_{m}\mathbf{\Sigma \Pi }_{m}\right\Vert ^{2}
&\leq &C_{2}\sum\limits_{i=1}^{n}\sum\limits_{i^{\prime }=1}^{n}\left[
\mathbb{E}\left( X\left( t_{i}\right) -\widetilde{X}\left( t_{i}\right)
\right) ^{2}+\mathbb{E}\left( X\left( t_{i^{\prime }}\right) -\widetilde{X}%
\left( t_{i^{\prime }}\right) \right) ^{2}\right] \\
&\leq &2C_{2}n^{2}\frac{1}{n}\sum\limits_{i=1}^{n}\mathbb{E}\left( X\left(
t_{i}\right) -\widetilde{X}\left( t_{i}\right) \right) ^{2}.
\end{eqnarray*}%
Since $n$ is fixed and the asymptotics are given with respect to $N$, the
number of replications of the process, the rate of convergence relies on the
quadratic error of the expansion of the process.

 For example take $d=1$, $T = [a,b]$,  $\mathcal{M} = \mathcal{M}_{N} = \left\{ (m) = \{1,\ldots,m\}, m=1,\ldots,N \right\}$, and for a process $X\left( t\right) $ with $t\in \left[a,b]%
\right] $,  consider its {Karhunen-Lo\`{e}ve expansion} (see for instance
\cite{adler}), i.e. write
$$
X\left( t\right) =\sum\limits_{\lambda =1}^{\infty
}Z_{\lambda}g_{\lambda}\left( t\right),
$$
where $Z_{\lambda}$ are centered random
variables with $\mathbb{E}\left( Z_{\lambda}^{2}\right)
=\gamma_{\lambda}^{2}$, where $\gamma_{\lambda}^{2}$ is the eigenvalue corresponding
to the eigenfunction $g_{\lambda}$ of the operator $\left( K f\right) \left(
t\right) =\int\limits_{a}^{b}\sigma \left( s,t\right) f\left( s\right) ds.$
If $X\left( t\right) $ is a Gaussian process then the random variables $%
Z_{\lambda} $ are Gaussian and stochastically independent. Hence, a natural approximation of $X\left( t\right) $ is given by
\begin{equation*}
\widetilde{X}\left( t\right) =\sum\limits_{\lambda =1}^{m}Z_{\lambda}g_{\lambda}\left( t\right)
.
\end{equation*}%
So we have that
\begin{equation*}
\mathbb{E}\left( X\left( t\right) -\widetilde{X}\left( t\right) \right) ^{2}=%
\mathbb{E}\left( \sum\limits_{\lambda =m+1}^{\infty }Z_{\lambda}g_{\lambda}\left( t\right)
\right) ^{2}=\sum\limits_{\lambda =m+1}^{\infty }\gamma_{\lambda}^{2}g_{\lambda}^{2}\left(
t\right) .
\end{equation*}%
therefore, if $\left\Vert g_{\lambda}\right\Vert _{L_{2}([a,b])}^{2}=1$ then $\mathbb{E}%
\left\Vert X\left( t\right) -\widetilde{X}\left( t\right) \right\Vert
_{L_{2}([a,b])}^{2}=\sum\limits_{l=m+1}^{\infty }\gamma_{\lambda}^{2}.$ Assume that the $\gamma_{\lambda}$'s have a polynomial decay of rate $\alpha >0$, namely $\gamma_{\lambda} \sim \lambda^{-\alpha}$, then we get an approximation error of
order $O\left( \left( m+1\right) ^{-2\alpha }\right) .$ Hence, we get that (under appropriate conditions on the design points $t_1,\ldots,t_n$)
\begin{equation*}
\left\Vert \mathbf{\Sigma -\Pi }_{m}\mathbf{\Sigma \Pi }_{m}\right\Vert
^{2}=O\left( \left( m+1\right) ^{-2\alpha }\right) .
\end{equation*}%
Finally, since in this example $\mathbb{E}\left\Vert \mathbf{\Sigma }-%
\widetilde{\mathbf{\Sigma }}\right\Vert ^{2}\leq K\left( \theta \right)
\underset{m\in \mathcal{M}_{N}}{\inf }\left( \left\Vert \mathbf{\Sigma -\Pi }%
_{m}\mathbf{\Sigma \Pi }_{m}\right\Vert ^{2}+\frac{\delta _{m}^{2}m}{N}%
\right) +O\left( \frac{1}{N}\right) $ then the quadratic risk is of order $%
N^{-\frac{2\alpha }{2\alpha +1}}$ as soon as $m\sim N^{1/(2\alpha +1)}$
belongs to the collection of models $\mathcal{M}_{N}$. In another framework,
if we consider a spline expansion, the rate of convergence for the
approximation given in \cite{Elogne} are of the same order.

Hence
we have obtained a model selection procedure which enables to recover the
best covariance model among a given collection. This method works without
strong assumptions on the process, in particular stationarity is not
assumed, but at the expand of necessary i.i.d observations of the process at
the same points. However the range of applications in broad, especially in
geophysics or epidemiology.

\section{Model Selection for Multidimensional Regression} \label{s:3}

\subsection{Oracle Inequality for multidimensional regression model}

\label{s:auxBaraud} Recall that we consider the following model
\begin{equation*}
\mathbf{x}_{i}\mathbf{x}_{i}^{T}=\mathbf{G\Psi G}^{T}+\mathbf{U}_{i},\quad
i=1,...,N,
\end{equation*}%
with i.i.d. matrix errors $\mathbf{U}_{i}$, $\mathbb{E}\left( \mathbf{U}%
_{i}\right) =\mathbf{0}$. This model can be equivalently rewritten in
vectorized form in the following way
\begin{equation*}
\mathbf{y=A\beta +u,}
\end{equation*}%
where $\mathbf{y}$ is a data vector, $\mathbb{E}\left( \mathbf{u}\right) =%
\mathbf{0}$, $\mathbf{A}$ is a known fixed matrix, and $\mathbf{\beta =}%
vech\left( \mathbf{\Psi }\right) $ is an unknown vector parameter. It is
worth of noting that this regression model has several peculiarities in
comparison with standard ones.\newline
$i)$ The error $\mathbf{u}$ has a specific correlation structure, namely $%
\mathbf{I}_{N}\otimes \mathbf{\Phi ,}$ where $\mathbf{\Phi }=V\left(
vec\left( \mathbf{x}_{i}\mathbf{x}_{i}^{T}\right) \right) $.\newline
$ii)$ In contrast with standard multivariate models, each coordinate of $%
\mathbf{y}$ depends on all the coordinates of $\mathbf{\beta }$.\newline
$iii)$ For any estimator $\widehat{\mathbf{\Sigma }}=\mathbf{G}\widehat{%
\mathbf{\Psi }}\mathbf{G}^{T}$ that be a linear function of the sample
covariance $\mathbf{S}$ of the data $\mathbf{x}_{1}$,...,$\mathbf{x}_{N}$
(and so, in particular, for the estimator minimizing $L_{N}$) it is possible
to construct an unbiased estimator of its quadratic risk $\mathbb{E}%
\left\Vert \mathbf{\Sigma -}\widehat{\mathbf{\Sigma }}\right\Vert ^{2}$.%

Assume we observe $\mathbf{y}_i$, $i=1,\dots,N$ random vectors of $\bR^d$
such that
\begin{equation}  \label{vectorizedmodel}
\mathbf{y}_{i}=\mathbf{f}^{i}\mathbf{+\varepsilon}_{i},\quad i=1,...,N,
\end{equation}
where $\mathbf{f}^{i}\mathbf{\in}\mathbb{R}^{d}$ are nonrandom and $\mathbf{%
\varepsilon}_{1},...,\mathbf{\varepsilon}_{N}$ are i.i.d. random vectors in $%
\mathbb{R}^{d}$ with $E\left( \mathbf{\varepsilon}_{1}\right) =\mathbf{0}$
and ${V} \left( \mathbf{\varepsilon}_{1}\right) =\mathbf{\Phi}$. For sake of
simplicity, we identify the function $g:\mathcal{X}\rightarrow\mathbb{R}^{d}$
with vectors $\left(g\left( x_{1}\right) \dots g\left(x_{N}\right)\right)^T
\in\mathbb{R}^{Nd}$ and we denote by $\left\langle a,b\right\rangle _{N}=%
\frac{1}{N}\sum\limits_{i=1}^{N}a_{i}^{T}b_{i}$, with $a=\left(a_{1}\dots
a_{N}\right)^T $ and $a_{i}\in\mathbb{R}^{d}$, the inner product of $\mathbb{%
R}^{Nd}$ associated to the norm $\left\Vert .\right\Vert _{N}$.

Given $N,d\in \mathbb{N}$, let $\left( \mathcal{L}_{m}\right) _{(m)\in
\mathcal{M}}$ be a finite family of linear subspaces of $\mathbb{R}^{Nd}$.
For each $(m)\in \mathcal{M}$, assume $\mathcal{L}_{m}$ has dimension $%
D_{m}\geq 1$. For each $(m)\in \mathcal{M}$, let $\widehat{\mathbf{f}}_{m}$
be the least squares estimator of $\mathbf{f=}\left( \left( \mathbf{f}%
^{1}\right) ^{T},...,\left( \mathbf{f}^{N}\right) ^{T}\right) ^{T}$ based on
the data $\mathbf{y=}\left( \mathbf{y}_{1},...,\mathbf{y}_{N}\right) $ under
the model $\mathcal{L}_{m}$; i.e.,
\begin{equation*}
\widehat{\mathbf{f}}_{m}=\arg \underset{\mathbf{v\in }\mathcal{L}_{m}}{\min }%
\left\{ \left\Vert \mathbf{y-v}\right\Vert _{N}^{2}\right\} =\mathbf{P}_{m}%
\mathbf{y,}
\end{equation*}%
where $\mathbf{P}_{m}$ is the projector matrix from $\mathbb{R}^{Nd}$ on $%
\mathcal{L}_{m}$. Write
\begin{align*}
\delta _{m}^{2}& =\frac{\mathrm{Tr}\left( \mathbf{P}_{m}\left( \mathbf{I}%
_{N}\otimes \mathbf{\Phi }\right) \right) }{D_{m}}, \\
\delta _{\sup }^{2}& =\max \left\{ \delta _{m}^{2}:m\in \mathcal{M}\right\} .
\end{align*}%
Given $\theta >0$, define the penalized estimator $\widetilde{\mathbf{f}}=%
\widehat{\mathbf{f}}_{\widehat{m}}$ , where%
\begin{equation*}
\widehat{m}=\arg \underset{(m)\in \mathcal{M}}{\min }\left\{ \left\Vert
\mathbf{y-}\widehat{\mathbf{f}}_{m}\right\Vert _{N}^{2}+pen\left( m\right)
\right\} ,
\end{equation*}%
with%
\begin{equation*}
pen\left( m\right) =\left( 1+\theta \right) \frac{\delta _{m}^{2}D_{m}}{N}.
\end{equation*}%
\begin{prop}
\label{ExtTh3.1Baraud}:
Let $q>0$ be given such that there exists $p>2\left(  1+q\right)  $ satisfying
$\mathbb{E}\left\Vert \mathbf{\varepsilon}_{1}\right\Vert ^{p}<\infty$. Then,
for some constants $K\left(  \theta\right)  >1$ and $c\left(  \theta
,p,q\right)  >0$ we have that%
\begin{equation}
\mathbb{E}\left(  \left\Vert \mathbf{f}-\widetilde{\mathbf{f}%
}\right\Vert_N^{2}-K\left(  \theta\right)  \mathcal{M}^{\ast}\right)  _{+}^{q}%
\leq\Delta_{p}^{q}\frac{\delta_{\sup}^{2q}}{N^{q}},\label{BoundForE}%
\end{equation}
where%
\begin{align*}
\Delta_{p}^{q}  & =C\left(  \theta,p,q\right)  \mathbb{E}\left\Vert
\mathbf{\varepsilon}_{1}\right\Vert ^{p}\left(  \sum\limits_{m\in\mathcal{M}%
}\delta_{m}^{-p}D_{m}^{-\left(  p/2-1-q\right)  }\right)  ,\\
\mathcal{M}^{\ast}  & =\underset{(m)\in\mathcal{M}}{\inf}\left\{  \left\Vert \mathbf{f-P}_{m}\mathbf{f}\right\Vert_N^{2}+\frac{\delta_{m}%
^{2}D_{m}}{N}\right\}  .
\end{align*}
\end{prop}This theorem is equivalent to Theorem~\ref{TeorConcentIneqLSEst}
using the vectorized version of the model \eqref{vectorizedmodel} and turns
to be an extension of Theorem 3.1 in \cite{Baraud00} to the multivariate
case. In a similar way, the following result constitutes also a natural
extension of Corollary 3.1 in \cite{Baraud00}. It is also closely related to the recent work in \cite{MR2471290}.
\begin{cor}
\label{ExtCor3.1Baraud}. Under
the assumptions of Proposition~\ref{ExtTh3.1Baraud} it holds that
\[
\left(  \mathbb{E}  \left\Vert \mathbf{f}-\widetilde
{\mathbf{f}}\right\Vert_N^{2q}\right)  ^{1/q}\leq2^{\left(
q^{-1}-1\right)  _{+}}\left[  K\left(  \theta\right)  \inf_{m\in\mathcal{M}%
}\left(  \left\Vert \mathbf{f-P}_{m}\mathbf{f}\right\Vert ^{2}+\frac
{\delta_{m}^{2}D_{m}}{N}\right)  +\frac{\Delta_{p}}{N}\delta_{\sup}%
^{2}\right]  ,
\]
where $\Delta_{p}$ was defined in Proposition~\eqref{ExtTh3.1Baraud}.
\end{cor}Under regularity assumptions for the function $\mathbf{f}$,
depending on a smoothness parameter $s$, the bias term is of order
\begin{equation*}
\left\Vert \mathbf{f-P}_{m}\mathbf{f}\right\Vert ^{2}=O(D_{m}^{-2s}).
\end{equation*}%
Hence, for $q=1$ we obtain the usual rate of convergence $N^{-\frac{2s}{2s+1}%
}$ for the quadratic risk as soon as the optimal choice $D_{m}=N^{\frac{1}{%
2s+1}}$ belongs to the collection of models, yielding the optimal rate of
convergence for the penalized estimator.

\subsection{Concentration Bound for multidimensional random process}

\label{s:concentration} These results are $d$-variate extensions of results
in \cite{Baraud00} (which are covered when $d=1$). Their proofs are deferred
to the Appendix.

\begin{prop}
\label{ExtCor5.1Baraud}(Extension of Corollary 5.1 in \cite{Baraud00}). Given
$N,d\in\mathbb{N}$, let $\widetilde{\mathbf{A}}\in\mathbb{R}^{Nd\times
Nd}\diagdown\left\{  \mathbf{0}\right\}  $ be a n.n.d. matrix and
$\mathbf{\varepsilon}_{1},...,\mathbf{\varepsilon}_{N}$ \ i.i.d random vectors
in $\mathbb{R}^{d}$ with $\mathbb{E}\left(  \mathbf{\varepsilon}_{1}\right)
=0$ and ${V}\left(  \mathbf{\varepsilon}_{1}\right)  =\mathbf{\Phi}$.
Write $\mathbf{\varepsilon}=\left(  \mathbf{\varepsilon}_{1}^{T}%
,...,\mathbf{\varepsilon}_{N}^{T}\right)  ^{T}$, $\zeta\left(
\mathbf{\varepsilon}\right)  =\sqrt{\mathbf{\varepsilon}^{T}\widetilde
{A}\mathbf{\varepsilon}}$, and $\gamma^{2}={\rm Tr}\left(  \widetilde{\mathbf{A}%
}\left(  \mathbf{I}_{N}\otimes\mathbf{\Phi}\right)  \right)  =\delta
^{2}{\rm Tr}\left(  \widetilde{\mathbf{A}}\right)  $. For all $p\geq2$ such that
$\mathbb{E}\left\Vert \mathbf{\varepsilon}_{1}\right\Vert ^{p}<\infty$ it
holds that, for all $x>0$
\begin{equation}
\mathbb{P}\left(  \zeta^{2}\left(  \mathbf{\varepsilon}\right)  \geq\delta
^{2}{\rm Tr}\left(  \widetilde{\mathbf{A}}\right)  +2\delta^{2}\sqrt{{\rm Tr}\left(
\widetilde{\mathbf{A}}\right)  \delta x}+\delta^{2}{\rm Tr}\left(  \widetilde
{\mathbf{A}}\right)  x\right)  \leq C\left(  p\right)  \frac{\mathbb{E}%
\left\Vert \varepsilon_{1}\right\Vert ^{p}{\rm Tr}\left(  \widetilde{\mathbf{A}%
}\right)  }{\delta^{p}\rho\left(  \widetilde{\mathbf{A}}\right)  x^{p/2}%
},\label{BoundCor5.1Baraud}%
\end{equation}
where the constant $C\left(  p\right)  $ depends only on $p.$
\end{prop}

Proposition~\ref{ExtCor5.1Baraud} reduces to Corollary 5.1 in \cite{Baraud00}
when when we only consider $d=1$, in which case $\delta^{2}=\left( \mathbf{%
\Phi }\right) _{11}=\sigma^{2}$ is the variance of the univariate i.i.d.
errors $\mathbf{\varepsilon}_{i}.$ 

\section{Appendix} \label{s:appen}

\subsection{Proofs of Preliminar results}

Proof of Proposition~\ref{pro:preliminar}
\begin{proof}
$a)$ The minimization problem posed in this theorem is equivalent to minimize
\[
h\left(  \mathbf{\Psi}\right)  =\left\Vert \overline{\mathbf{Y}}-\mathbf{G\Psi
G}^{T}\right\Vert ^{2}.
\]

The Frobenius norm $\left\Vert {.}\right\Vert $ is invariant by the $vec$
operation. Furthermore,$\mathbf{\Psi\in}\mathcal{S}_{m}$ can be represented by
means of $\mathbf{\delta=}vec\left(  \mathbf{\Psi}\right)  =\mathbf{D}%
_{q}\mathbf{\beta}$ where $\mathbf{\beta\in}\mathbb{R}^{q\left(  q+1\right)
/2}$. These facts and the identity%
\begin{equation}
vec\left(  \mathbf{ABC}\right)  =\left(  \mathbf{C}^{T}\otimes\mathbf{A}%
\right)  vec\left(  \mathbf{B}\right)  \label{vec(ABC)}%
\end{equation}
allow one to rewrite
\[
h\left(  \mathbf{\Psi}\right)  =\left\Vert \overline{\mathbf{y}}-\left(
\mathbf{G}\otimes\mathbf{G}\right)  \mathbf{D}_{q}\mathbf{\beta}\right\Vert
^{2},
\]
where $\overline{\mathbf{y}}=vec\left(  \overline{\mathbf{Y}}\right)  $.
Minimization of this quadratic function with respect to $\mathbf{\beta}$
in\textbf{\ }$\mathbb{R}^{q\left(  q+1\right)  /2}$ is equivalent to solve the
normal equation%
\[
\mathbf{D}_{q}^{T}\left(  \mathbf{G}\otimes\mathbf{G}\right)  ^{T}\left(
\mathbf{G}\otimes\mathbf{G}\right)  \mathbf{D}_{q}\mathbf{\beta=D}_{q}%
^{T}\left(  \mathbf{G}\otimes\mathbf{G}\right)  ^{T}\overline{\mathbf{y}}.
\]
By using the identities
\[
\mathbf{D}_{q}^{T}vec\left(  \mathbf{A}\right)  =vech\left(  \mathbf{A+A}%
^{T}-diag\left(  \mathbf{A}\right)  \right)
\]
and \ref{vec(ABC)}, said normal equation can be rewritten
\[
vech\left(  \mathbf{G}^{T}\mathbf{G}\left(  \mathbf{\Psi+\Psi}^{T}\right)
\mathbf{G}^{T}\mathbf{G-}diag\left(  \mathbf{G}^{T}\mathbf{G\Psi G}%
^{T}\mathbf{G}\right)  \right)  =vech\left(  \mathbf{G}^{T}\left(
\overline{\mathbf{Y}}+\overline{\mathbf{Y}}^{T}\right)  \mathbf{G}\right)  .
\]
Finally, it can be verified that $\widehat{\mathbf{\Psi}}$ given by
(\ref{GeneralPhiTecho}) satisfies this equation as a consequence of the fact
that such $\widehat{\mathbf{\Psi}}$ it holds that%
\[
\mathbf{G}^{T}\mathbf{G\widehat{\mathbf{\Psi}}G}^{T}\mathbf{G}=vech\left(
\mathbf{G}^{T}\left(  \frac{\overline{\mathbf{Y}}+\overline{\mathbf{Y}}^{T}%
}{2}\right)  \mathbf{G}\right)  .
\]
$b)$ It straightforwardly follows from part $a)$.
\end{proof}

\subsection{Proofs of Main Results}

Proof of Proposition~\eqref{ExtTh3.1Baraud}
\begin{proof}
The proof follows the guidelines of the proof in \cite{Baraud00}. More generally we will prove that for any $\eta>0$ and any sequence of
positive numbers $L_{m}$, if the penalty function $pen:$ $\mathcal{M}%
\longrightarrow\mathbb{R}_{+}$ is chosen to satisfy:%
\begin{equation}
pen\left(  m\right)  =\left(  1+\eta+L_{m}\right)  \frac{\delta_{m}^{2}}%
{N}D_{m}\text{ for all }(m)\in\mathcal{M},\label{pen}%
\end{equation}
then for each $x>0$ and $p\geq2$%
\begin{equation}
\mathbb{P}\left(  \mathcal{H}\left(  \mathbf{f}\right)  \geq\left(  1+\frac
{2}{\eta}\right)  \frac{x}{N}\delta_{m}^{2}\right)  \leq c\left(
p,\eta\right)  \mathbb{E}\left\Vert \mathbf{\varepsilon}_{1}\right\Vert
^{p}\sum\limits_{(m)\in\mathcal{M}}\frac{1}{\delta_{m}^{p}}\frac{D_{m}\vee
1}{\left(  L_{m}D_{m}+x\right)  ^{p/2}},\label{ineqPH}%
\end{equation}
where we have set
\[
\mathcal{H}\left(  \mathbf{f}\right)  =\left[  \left\Vert \mathbf{f}%
-\widetilde{\mathbf{f}}\right\Vert _{N}^{2}-\left(  2-\frac{4}{\eta}\right)
\underset{(m)\in\mathcal{M}}{\inf}\left\{  d_{N}^{2}\left(  \mathbf{f}%
,\mathcal{L}_{m}\right)  +pen\left(  m\right)  \right\}  \right]  _{+}.
\]

To obtain (\ref{BoundForE}), take $\eta=\frac{\theta}{2}=L_{m}$. As for each
$(m)\in\mathcal{M}$,
\begin{align*}
d_{N}^{2}\left(  \mathbf{f},\mathcal{L}_{m}\right)  +pen\left(  m\right)   &
\leq d_{N}^{2}\left(  \mathbf{f},\mathcal{L}_{m}\right)  +\left(
1+\theta\right)  \frac{\delta_{m}^{2}}{N}D_{m}\\
& \leq\left(  1+\theta\right)  \left(  d_{N}^{2}\left(  \mathbf{f}%
,\mathcal{L}_{m}\right)  +\frac{\delta_{m}^{2}}{N}D_{m}\right)
\end{align*}
we get that for all $q>0$,%
\begin{equation}
\mathcal{H}^{q}\left(  \mathbf{f}\right)  \geq\left[  \left\Vert
\mathbf{f}-\widetilde{\mathbf{f}}\right\Vert _{N}^{2}-\left(  2+\frac
{8}{\theta}\right)  \left(  1+\theta\right)  \mathcal{M}^{\ast}\right]  _{+}%
^{q}=\left[  \left\Vert \mathbf{f}-\widetilde{\mathbf{f}}\right\Vert _{N}%
^{2}-K\left(  \theta\right)  \mathcal{M}^{\ast}\right]  _{+}^{q},\label{ineqHq}%
\end{equation}
where $K\left(  \theta\right)  =\left(  2+\frac{8}{\theta}\right)  \left(
1+\theta\right)  $.

Since%
\[
\mathbb{E}\left(  \mathcal{H}^{q}\left(  \mathbf{f}\right)  \right)
=\int\limits_{0}^{\infty}qu^{q-1}\mathbb{P}\left(  \mathcal{H}\left(
\mathbf{f}\right)  >u\right)  du,
\]
we derive from (\ref{ineqHq}) and (\ref{ineqPH}) that for all $p>2\left(
1+q\right)  $%
\begin{align*}
\mathbb{E}\left[  \left(  \left\Vert \mathbf{f}-\widetilde{\mathbf{f}%
}\right\Vert _{N}^{2}-K\left(  \theta\right)  \mathcal{M}^{\ast}\right)  _{+}%
^{q}\right]   & \leq\mathbb{E}\left(  \mathcal{H}^{q}\left(  \mathbf{f}%
\right)  \right) \\
& \leq c\left(  p,\theta\right)  \left(  1+\frac{4}{\theta}\right)  ^{q}%
\frac{\mathbb{E}\left\Vert \mathbf{\varepsilon}_{1}\right\Vert ^{p}}{N^{q}%
}\sum\limits_{m\in\mathcal{M}}\frac{\delta_{m}^{2q}}{\delta_{m}^{p}}%
\int\limits_{0}^{\infty}qx^{q-1}\left[  \frac{D_{m}\vee1}{\left(  \frac
{\theta}{2}D_{m}+x\right)  ^{p/2}}\wedge1\right]  dx\\
& \leq c^{\prime}\left(  p,q,\theta\right)  \frac{\mathbb{E}\left\Vert
\mathbf{\varepsilon}_{1}\right\Vert ^{p}}{N^{q}}\delta_{\sup}^{2q}\left[
\sum\limits_{(m)\in\mathcal{M}}\delta_{m}^{-p}D_{m}^{-\left(  p/2-1-q\right)
}\right]
\end{align*}
using that $\mathbb{P}\left(  \mathcal{H}\left(  \mathbf{f}\right)  >u\right)
\leq1$.\\ Indeed, for $m\in\mathcal{M}$ such that $D_{m}\geq1$, using that $q-1-p/2<0 $,
we get the following bounds 
\begin{align}
\frac{\delta_{m}^{2q}}{\delta_{m}^{p}}\int\limits_{0}^{\infty}qx^{q-1}\left[
\frac{D_{m}\vee1}{\left(  \frac{\theta}{2}D_{m}+x\right)  ^{p/2}}%
\wedge1\right]  dx  & \leq\delta_{\sup}^{2q}\delta_{m}^{-p}\int\limits_{0}%
^{\infty}qx^{q-1}\left[  \frac{D_{m}}{\left(  \frac{\theta}{2}D_{m}+x\right)
^{p/2}}\right]  dx\nonumber\\
 =\delta_{\sup}^{2q}\delta_{m}^{-p} & \left(  \int\limits_{0}^{D_{m}}qx^{q-1}
\left[  \frac{D_{m}}{\left(  \frac{\theta}{2}D_{m}+x\right)  ^{p/2}}\right]
dx+\int\limits_{D_{m}}^{\infty}qx^{q-1}\left[  \frac{D_{m}}{\left(
\frac{\theta}{2}D_{m}+x\right)  ^{p/2}}\right]  dx\right) \nonumber\\
& \leq\delta_{\sup}^{2q}\delta_{m}^{-p}\left(  \frac{D_{m}}{\left(
\frac{\theta}{2}D_{m}\right)  ^{p/2}}\int\limits_{0}^{D_{m}}qx^{q-1}%
dx+D_{m}\int\limits_{D_{m}}^{\infty}qx^{q-1}\left[  \frac{1}{x^{p/2}}\right]
dx\right) \nonumber\\
& =\delta_{\sup}^{2q}\delta_{m}^{-p}\left(  2^{p/2}\theta^{-p/2}D_{m}%
^{1-p/2}\int\limits_{0}^{D_{m}}qx^{q-1}dx+D_{m}\int\limits_{D_{m}}^{\infty
}qx^{q-1-p/2}dx\right) \nonumber\\
& =\delta_{\sup}^{2q}\delta_{m}^{-p}\left(  2^{p/2}\theta^{-p/2}D_{m}%
^{1-p/2}\left[  D_{m}^{q}\right]  +D_{m}\left[  \frac{q}{p/2-q}D_{m}%
^{q-p/2}\right]  \right) \nonumber\\
& =\delta_{\sup}^{2q}\delta_{m}^{-p}\left(  2^{p/2}\theta^{-p/2}%
D_{m}^{1-p/2+q}+D_{m}^{1-p/2+q}\left[  \frac{q}{p/2-q}\right]  \right)
\nonumber\\
& =\delta_{\sup}^{2q}\delta_{m}^{-p}\left(  D_{m}^{-\left(  p/2-1-q\right)  }
\left[  2^{p/2}\theta^{-p/2}+\frac{q}{p/2-q}\right]  \right)  .\label{Eq6}%
\end{align}

\eqref{Eq6} enables to conclude that \eqref{BoundForE} holds assuming \eqref{ineqPH}.\vskip .1in

We now turn to the proof of (\ref{ineqPH}). Recall that, we
identify the function $g:\mathcal{X}\rightarrow\mathbb{R}^{d}$ with vectors
$\left(g\left(  x_{1}\right) \dots g\left(x_{N}\right)\right)^T  \in\mathbb{R}^{Nd}$ and we define the empirical scalar product as $\left\langle a,b\right\rangle
_{N}=\frac{1}{N}\sum\limits_{i=1}^{N}a_{i}^{T}b_{i}$, with $a=\left(a_{1}\dots a_{N}\right)^T  $ and $a_{i}\in\mathbb{R}^{d}$, the inner product of $\mathbb{R}%
^{Nd}$ associated to the norm $\left\Vert .\right\Vert _{N}$. For each $(m)\in
\mathcal{M}$ we denote by $\mathbf{P}_{m}$ the orthogonal projector onto the linear
space $\left\{  \left(g\left(  x_{1}\right) \dots
g\left(  x_{N}\right)
\right)^T  :g\in\mathcal{L}_{m}\right\}  \subset$ $\mathbb{R}^{Nd}$. This linear
space is also denoted by $\mathcal{L}_{m}$. From now on, the subscript $m$
denotes any minimizer of the function $m^{\prime}\rightarrow$ $\left\Vert
\mathbf{f}-\mathbf{P}_{m^{\prime}}\mathbf{f}\right\Vert ^{2}+pen\left(
m^{\prime}\right)  $, $(m^{\prime})\in\mathcal{M}_{N}$. For any $\mathbf{g}\in$
$\mathbb{R}^{Nd}$ we define the least-squares loss function by
\[
\gamma_{N}\left(  \mathbf{g}\right)  =\left\Vert \mathbf{y-g}\right\Vert
_{N}^{2}%
\]
Using the definition of $\gamma_{N}$ we have that for all $\mathbf{g}\in$
$\mathbb{R}^{Nd}$,%
\[
\gamma_{N}\left(  \mathbf{g}\right)  =\left\Vert \mathbf{f}%
+\mathbf{\ \mathbf{\varepsilon}}-\mathbf{g}\right\Vert _{N}^{2}.
\]

Then we derive that%
\[
\left\Vert \mathbf{f}-\mathbf{g}\right\Vert _{N}^{2}=\gamma_{N}\left(
\mathbf{f}\right)  +2\left\langle \mathbf{f}-\mathbf{y}%
,\mathbf{\ \mathbf{\varepsilon}}\right\rangle _{N}+\left\Vert
\mathbf{\ \mathbf{\varepsilon}}\right\Vert _{N}^{2}%
\]
and therefore%
\begin{equation}
\left\Vert \mathbf{f}-\widetilde{\mathbf{f}}\right\Vert _{N}^{2}-\left\Vert
\mathbf{f}-\mathbf{P}_{m}\mathbf{f}\right\Vert _{N}^{2}=\gamma_{N}\left(
\widetilde{\mathbf{f}}\right)  -\gamma_{N}\left(  \mathbf{P}_{m}%
\mathbf{f}\right)  +2\left\langle \widetilde{\mathbf{f}}-\mathbf{P}%
_{m}\mathbf{f},\mathbf{\ \mathbf{\varepsilon}}\right\rangle _{N}.\label{eq7}%
\end{equation}

By the definition of $\widetilde{\mathbf{f}}$, we know that%
\[
\gamma_{N}\left(  \widetilde{\mathbf{f}}\right)  +pen\left(  \widehat
{m}\right)  \leq\gamma_{N}\left(  \mathbf{g}\right)  +pen\left(  m\right)
\]
for all $(m)\in\mathcal{M}$ and for all $\mathbf{g}\in\mathcal{L}_{m}$. Then%
\begin{equation}
\gamma_{N}\left(  \widetilde{\mathbf{f}}\right)  -\gamma_{N}\left(
\mathbf{P}_{m}\mathbf{f}\right)  \leq pen\left(  m\right)  -pen\left(
\widehat{m}\right)  .\label{eq7.1}%
\end{equation}
So we get from (\ref{eq7}) and (\ref{eq7.1}) that%
\begin{equation}
\left\Vert \mathbf{f}-\widetilde{\mathbf{f}}\right\Vert _{N}^{2}\leq\left\Vert
\mathbf{f}-\mathbf{P}_{m}\mathbf{f}\right\Vert _{N}^{2}+pen\left(  m\right)
-pen\left(  \widehat{m}\right)  +2\left\langle \mathbf{f}-\mathbf{P}%
_{m}\mathbf{f},\mathbf{\ \mathbf{\varepsilon}}\right\rangle _{N}+2\left\langle
\mathbf{P}_{\widehat{m}}\mathbf{f}-\mathbf{f},\mathbf{\ \mathbf{\varepsilon}%
}\right\rangle _{N}+2\left\langle \widetilde{\mathbf{f}}-\mathbf{P}%
_{\widehat{m}}\mathbf{f},\mathbf{\ \mathbf{\varepsilon}}\right\rangle
_{N}.\label{eq8}%
\end{equation}

In the following we set for each $(m^{\prime})\in\mathcal{M}$,
\begin{align*}
\mathcal{B}_{m^{\prime}}  & =\left\{  \mathbf{g}\in\mathcal{L}_{m^{\prime}%
}:\left\Vert \mathbf{g}\right\Vert _{N}\leq1\right\}  ,\\
G_{m^{\prime}}  & =\underset{t\in\mathcal{B}_{m^{\prime}}}{\sup}\left\langle
\mathbf{g},\mathbf{\ \mathbf{\varepsilon}}\right\rangle _{N}=\left\Vert
\mathbf{P}_{m^{\prime}}\mathbf{\ \mathbf{\varepsilon}}\right\Vert _{N},\\
& \mathbf{u}_{m^{\prime}}=\begin{cases}
\frac{\mathbf{P}_{m^{\prime}}\mathbf{f}-\mathbf{f}}{\left\Vert \mathbf{P}%
_{m^{\prime}}\mathbf{f}-\mathbf{f}\right\Vert _{N}} & \text{ if }\left\Vert
\mathbf{P}_{m^{\prime}}\mathbf{f}-\mathbf{f}\right\Vert _{N}\neq0\\
0 & \text{ otherwise.}
\end{cases}
\end{align*}

Since $\widetilde{\mathbf{f}}=$ $\mathbf{P}_{\widehat{m}}$ $\mathbf{f}+$
$\mathbf{P}_{\widehat{m}}$\textbf{\ $\mathbf{\varepsilon}$}, (\ref{eq8}) gives%
\begin{equation*}
\left\Vert \mathbf{f}-\widetilde{\mathbf{f}}\right\Vert _{N}^{2}\leq\left\Vert
\mathbf{f}-\mathbf{P}_{m}\mathbf{f}\right\Vert _{N}^{2}+pen\left(  m\right)
-pen\left(  \widehat{m}\right)
\end{equation*}
\begin{equation}
+2\left\Vert \mathbf{f}-\mathbf{P}%
_{m}\mathbf{f}\right\Vert _{N}\left\vert \left\langle \mathbf{u}%
_{m},\mathbf{\ \mathbf{\varepsilon}}\right\rangle _{N}\right\vert +2\left\Vert
\mathbf{f}-\mathbf{P}_{\widehat{m}}\mathbf{f}\right\Vert _{N}\left\vert
\left\langle \mathbf{u}_{\widehat{m}},\mathbf{\ \mathbf{\varepsilon}%
}\right\rangle _{N}\right\vert +2G_{\widehat{m}}^{2}.\label{eq9}%
\end{equation}

Using repeatedly the following elementary inequality that holds for all
positive numbers $\alpha,x,z$%
\begin{equation}
2xz\leq\alpha x^{2}+\frac{1}{\alpha}z^{2}\label{eq101}%
\end{equation}
we get for any $m^{\prime}\in\mathcal{M}$%
\begin{equation}
2\left\Vert \mathbf{f}-\mathbf{P}_{m^{\prime}}\mathbf{f}\right\Vert \left\vert
\left\langle \mathbf{u}_{m^{\prime}},\mathbf{\ \mathbf{\varepsilon}%
}\right\rangle _{N}\right\vert \leq\alpha\left\Vert \mathbf{f}-\mathbf{P}%
_{m^{\prime}}\mathbf{f}\right\Vert ^{2}+\frac{1}{\alpha}\left\vert
\left\langle \mathbf{u}_{m^{\prime}},\mathbf{\ \mathbf{\varepsilon}%
}\right\rangle _{N}\right\vert ^{2}.\label{eq11}%
\end{equation}

By Pythagoras Theorem we have
\begin{align}
\left\Vert \mathbf{f}-\widetilde{\mathbf{f}}\right\Vert _{N}^{2}  &
=\left\Vert \mathbf{f}-\mathbf{P}_{\widehat{m}}\mathbf{f}\right\Vert
_{N}^{2}+\left\Vert \mathbf{P}_{\widehat{m}}\mathbf{f}-\widetilde{\mathbf{f}%
}\right\Vert _{N}^{2}\nonumber\\
& =\left\Vert \mathbf{f}-\mathbf{P}_{\widehat{m}}\mathbf{f}\right\Vert
_{N}^{2}+G_{\widehat{m}}^{2}.\label{eq12}%
\end{align}

We derive from (\ref{eq9}) and (\ref{eq11}) that for any $\alpha>0$:%
\[
\left\Vert \mathbf{f}-\widetilde{\mathbf{f}}\right\Vert _{N}^{2}\leq\left\Vert
\mathbf{f}-\mathbf{P}_{m}\mathbf{f}\right\Vert _{N}^{2}+\alpha\left\Vert
\mathbf{f}-\mathbf{P}_{m}\mathbf{f}\right\Vert _{N}^{2}+\frac{1}{\alpha
}\left\langle \mathbf{u}_{m},\mathbf{\ \mathbf{\varepsilon}}\right\rangle
_{N}^{2} \]
\[+\alpha\left\Vert \mathbf{f}-\mathbf{P}_{\widehat{m}}\mathbf{f}%
\right\Vert _{N}^{2}+\frac{1}{\alpha}\left\langle \mathbf{u}_{\widehat{m}%
},\mathbf{\ \mathbf{\varepsilon}}\right\rangle _{N}^{2}+2G_{\widehat{m}}%
^{2}+pen\left(  m\right)  -pen\left(  \widehat{m}\right)  .
\]

Now taking into account that by equation (\ref{eq12}) $\left\Vert
\mathbf{f}-\mathbf{P}_{\widehat{m}}\mathbf{f}\right\Vert _{N}^{2}=\left\Vert
\mathbf{f}-\widetilde{\mathbf{f}}\right\Vert _{N}^{2}-G_{\widehat{m}}^{2}$ the
above inequality is equivalent to:%
\begin{equation*}
\left(  1-\alpha\right)  \left\Vert \mathbf{f}-\widetilde{\mathbf{f}%
}\right\Vert _{N}^{2}\leq\left(  1+\alpha\right)  \left\Vert \mathbf{f}%
-\mathbf{P}_{m}\mathbf{f}\right\Vert _{N}^{2}+\frac{1}{\alpha}\left\langle
\mathbf{u}_{m},\mathbf{\varepsilon}\right\rangle _{N}^{2}
\end{equation*}
\begin{equation}
+\frac{1}{\alpha
}\left\langle \mathbf{u}_{\widehat{m}},\mathbf{\varepsilon}\right\rangle
_{N}^{2}+\left(  2-\alpha\right)  G_{\widehat{m}}^{2}+pen\left(  m\right)
-pen\left(  \widehat{m}\right)  .\label{eq13}%
\end{equation}

We choose $\alpha=\frac{2}{2+\eta}\in\left]  0,1\right[  $, but for sake of simplicity
we keep using the notation $\alpha$. Let $\widetilde{p}_{1}$ and
$\widetilde{p}_{2}$ be two functions depending on $\eta$ mapping $\mathcal{M}$
into $\mathbb{R}_{+}$. They will be specified later to satisfy%
\begin{equation}
pen\left(  m^{\prime}\right)  \geq\left(  2-\alpha\right)  \widetilde{p}%
_{1}\left(  m^{\prime}\right)  +\frac{1}{\alpha}\widetilde{p}_{2}\left(
m^{\prime}\right)  \text{ }\forall (m^{\prime})\in\mathcal{M}_{.}\label{eq14}%
\end{equation}

Since $\frac{1}{\alpha}\widetilde{p}_{2}\left(  m^{\prime}\right)  \leq
pen\left(  m^{\prime}\right)  $ and $1+\alpha\leq2$, we get from (\ref{eq13})
and (\ref{eq14}) that%
\begin{align}
\left(  1-\alpha\right)  \left\Vert \mathbf{f}-\widetilde{\mathbf{f}%
}\right\Vert _{N}^{2}  & \leq \left(  1+\alpha\right)  \left\Vert
\mathbf{f}-\mathbf{P}_{m}\mathbf{f}\right\Vert _{N}^{2}+pen\left(  m\right)
+\frac{1}{\alpha}\widetilde{p}_{2}\left(  m\right)  +\left(  2-\alpha\right)
\left(  G_{\widehat{m}}^{2}-\widetilde{p}_{1}\left(  \widehat{m}\right)
\right) \nonumber \\
& +\frac{1}{\alpha}\left(  \left\langle \mathbf{u}_{\widehat{m}%
},\mathbf{\varepsilon}\right\rangle _{N}^{2}-\widetilde{p}_{2}\left(
\widehat{m}\right)  \right)  +\frac{1}{\alpha}\left(  \left\langle
\mathbf{u}_{m},\mathbf{\varepsilon}\right\rangle _{N}^{2}-\widetilde{p}%
_{2}\left(  m\right)  \right) \nonumber\\
& \leq2\left(  \left\Vert \mathbf{f}-\mathbf{P}_{m}\mathbf{f}\right\Vert
_{N}^{2}+pen\left(  m\right)  \right)  +\left(  2-\alpha\right)  \left(
G_{\widehat{m}}^{2}-\widetilde{p}_{1}\left(  \widehat{m}\right)  \right)  \nonumber \\
& +\frac{1}{\alpha}\left(  \left\langle \mathbf{u}_{\widehat{m}}%
,\mathbf{\varepsilon}\right\rangle _{N}^{2}-\widetilde{p}_{2}\left(
\widehat{m}\right)  \right)  +\frac{1}{\alpha}\left(  \left\langle
\mathbf{u}_{m},\mathbf{\varepsilon}\right\rangle _{N}^{2}-\widetilde{p}%
_{2}\left(  m\right)  \right)  .\label{eq15}%
\end{align}

As $\frac{2}{1-\alpha}=2+\frac{4}{\eta}$ we obtain that%
\begin{align*}
\left(  1-\alpha\right)  \mathcal{H}\left(  \mathbf{f}\right)   & =\left\{
\left(  1-\alpha\right)  \left\Vert \mathbf{f}-\widetilde{\mathbf{f}%
}\right\Vert _{N}^{2}-\left(  1-\alpha\right)  \left(  2+\frac{4}{\eta
}\right)  \underset{m^{\prime}\in \mathcal{M}}{\inf}\left(  \left\Vert
\mathbf{f}-\mathbf{P}_{m^{\prime}}\mathbf{f}\right\Vert _{N}^{2}+pen\left(
m^{\prime}\right)  \right)  \right\}  _{+}\\
& =\left\{  \left(  1-\alpha\right)  \left\Vert \mathbf{f}-\widetilde
{\mathbf{f}}\right\Vert _{N}^{2}-2\left(  \left\Vert \mathbf{f}-\mathbf{P}%
_{m}\mathbf{f}\right\Vert _{N}^{2}+2pen\left(  m\right)  \right)  \right\}
_{+}\\
& \leq\left\{  \left(  2-\alpha\right)  \left(  G_{\widehat{m}}^{2}%
-\widetilde{p}_{1}\left(  \widehat{m}\right)  \right)  +\frac{1}{\alpha
}\left(  \left\langle \mathbf{u}_{\widehat{m}},\mathbf{\varepsilon
}\right\rangle _{N}^{2}-\widetilde{p}_{2}\left(  \widehat{m}\right)  \right)
+\frac{1}{\alpha}\left(  \left\langle \mathbf{u}_{m},\mathbf{\varepsilon
}\right\rangle _{N}^{2}-\widetilde{p}_{2}\left(  m\right)  \right)  \right\}
_{+}%
\end{align*}
using that $m$ minimizes the function $\left\Vert \mathbf{f}-\mathbf{P}%
_{m^{\prime}}\right\Vert ^{2}+pen\left(  m^{\prime}\right)  $ and (\ref{eq15}).

For any $x>0,$%
\begin{align}
\mathbb{P}\left(  \left(  1-\alpha\right)  \mathcal{H}\left(  \mathbf{f}%
\right)  \geq\frac{x\delta_{m}^{2}}{N}\right)   & \leq\mathbb{P}\left(
\exists m^{\prime}\in \mathcal{M}:\left(  2-\alpha\right)  \left(  G_{m^{\prime}%
}^{2}-\widetilde{p}_{1}\left(  m^{\prime}\right)  \right)  \geq\frac
{x\delta_{m^{\prime}}^{2}}{3N}\right) \nonumber\\
& +\mathbb{P}\left(  \exists m^{\prime}\in \mathcal{M}:\frac{1}{\alpha}\left(
\left\langle \mathbf{u}_{m^{\prime}},\mathbf{\varepsilon}\right\rangle
_{N}^{2}-\widetilde{p}_{2}\left(  m^{\prime}\right)  \right)  \geq
\frac{x\delta_{m^{\prime}}^{2}}{3N}\right) \nonumber\\
& \leq\sum\limits_{m^{\prime}\in \mathcal{M}}\mathbb{P}\left(  \left(  2-\alpha\right) \left( \left\Vert
\mathbf{P}_{m^{\prime}}\mathbf{\varepsilon}\right\Vert _{N}^{2}-\widetilde
{p}_{1}\left(  m^{\prime}\right) \right)  \geq\frac{x\delta_{m^{\prime}}^{2}}%
{3N}\right) \nonumber\\
& +\sum\limits_{m^{\prime}\in \mathcal{M}}\mathbb{P}\left(  \frac{1}{\alpha}\left(
\left\langle \mathbf{u}_{m^{\prime}},\mathbf{\varepsilon}\right\rangle
_{N}^{2}-\widetilde{p}_{2}\left(  m^{\prime}\right)  \right)  \geq
\frac{x\delta_{m^{\prime}}^{2}}{3N}\right) \nonumber\\
& :=\sum\limits_{m^{\prime}\in \mathcal{M}}P_{1,m^{\prime}}\left(  x\right)
+\sum\limits_{m^{\prime}\in \mathcal{M}}P_{2,m^{\prime}}\left(  x\right)
.\label{eq16}%
\end{align}

We first bound $P_{2,m^{\prime}}\left(  x\right)  $. Let $t$ be some positive
number,%
\begin{equation}
\mathbb{P}\left(  \left\vert \left\langle \mathbf{u}_{m^{\prime}%
},\mathbf{\varepsilon}\right\rangle _{N}\right\vert \geq t\right)  \leq
t^{-p}\mathbb{E}\left(  \left\vert \left\langle \mathbf{u}_{m^{\prime}%
},\mathbf{\varepsilon}\right\rangle _{N}\right\vert ^{p}\right)  .\label{eq17}%
\end{equation}

Since $\left\langle \mathbf{u}_{m^{\prime}},\mathbf{\varepsilon}\right\rangle
_{N}=\frac{1}{N}\sum\limits_{i=1}^{N}\left\langle \mathbf{u}_{im^{\prime}%
},\mathbf{\varepsilon}_{i}\right\rangle $ with $\mathbf{\varepsilon}_{i}$
i.i.d. and with zero mean, then by Rosenthal's inequality we know that for
some constant $c\left(  p\right)  $ that depends on $p$ only
\begin{align}
c^{-1}\left(  p\right)  N^{p}\mathbb{E}\left\vert \left\langle \mathbf{u}%
_{m^{\prime}},\mathbf{\varepsilon}\right\rangle _{N}\right\vert ^{p}  &
\leq\sum\limits_{i=1}^{N}\mathbb{E}\left\vert \left\langle \mathbf{u}%
_{im^{\prime}},\mathbf{\varepsilon}_{i}\right\rangle \right\vert ^{p}+\left(
\sum\limits_{i=1}^{N}\mathbb{E}\left(  \left\langle \mathbf{u}_{im^{\prime}%
},\mathbf{\varepsilon}_{i}\right\rangle ^{2}\right)  \right)  ^{\frac{p}{2}%
}\nonumber\\
& \leq\sum\limits_{i=1}^{N}\mathbb{E}\left\Vert \mathbf{u}_{im^{\prime}%
}\right\Vert ^{p}\left\Vert \mathbf{\varepsilon}_{i}\right\Vert ^{p}+\left(
\sum\limits_{i=1}^{N}\mathbb{E}\left\Vert \mathbf{u}_{im^{\prime}}\right\Vert
^{2}\left\Vert \mathbf{\varepsilon}_{i}\right\Vert ^{2}\right)  ^{\frac{p}{2}%
}\nonumber\\
& =\mathbb{E}\left\Vert \mathbf{\varepsilon}_{1}\right\Vert ^{p}%
\sum\limits_{i=1}^{N}\left\Vert \mathbf{u}_{im^{\prime}}\right\Vert
^{p}+\left(  \mathbb{E}\left\Vert \mathbf{\varepsilon}_{1}\right\Vert
^{2}\right)  ^{\frac{p}{2}}\left(  \sum\limits_{i=1}^{N}\left\Vert
\mathbf{u}_{im^{\prime}}\right\Vert ^{2}\right)  ^{\frac{p}{2}}.\label{eq18}%
\end{align}

Since $p\geq2$, $\left(  \mathbb{E}\left\Vert \mathbf{\varepsilon}%
_{1}\right\Vert ^{2}\right)  ^{\frac{1}{2}}\leq\left(  \mathbb{E}\left\Vert
\mathbf{\varepsilon}_{1}\right\Vert ^{p}\right)  ^{\frac{1}{p}}$ and
\begin{equation}
\left(  \mathbb{E}\left\Vert \mathbf{\varepsilon}_{1}\right\Vert ^{2}\right)
^{\frac{p}{2}}\leq\mathbb{E}\left\Vert \mathbf{\varepsilon}_{1}\right\Vert
^{p}.\label{eq19}%
\end{equation}

Using also that by definition $\left\Vert \mathbf{u}_{m^{\prime}}\right\Vert
_{N}^{2}=\frac{1}{N}\sum\limits_{i=1}^{N}\left\Vert \mathbf{u}_{im^{\prime}%
}\right\Vert ^{2}=1$, then $\frac{\left\Vert \mathbf{u}_{im^{\prime}%
}\right\Vert ^{2}}{N}\leq1$ and therefore $\frac{\left\Vert \mathbf{u}%
_{im^{\prime}}\right\Vert }{N^{\frac{1}{2}}}\leq1$. Thus%
\begin{equation}
\sum\limits_{i=1}^{N}\left\Vert \mathbf{u}_{im^{\prime}}\right\Vert
^{p}=N^{\frac{p}{2}}\sum\limits_{i=1}^{N}\left(  \frac{\left\Vert
\mathbf{u}_{im^{\prime}}\right\Vert }{N^{\frac{1}{2}}}\right)  ^{p}\leq
N^{\frac{p}{2}}\sum\limits_{i=1}^{N}\left(  \frac{\left\Vert \mathbf{u}%
_{im^{\prime}}\right\Vert }{N^{\frac{1}{2}}}\right)  ^{2}=N^{\frac{p}{2}%
}\left\Vert \mathbf{u}_{m^{\prime}}\right\Vert _{N}^{2}=N^{\frac{p}{2}%
}.\label{eq20}%
\end{equation}

We deduce from (\ref{eq18}), (\ref{eq19}) and (\ref{eq20}) that%
\[
c^{-1}\left(  p\right)  N^{p}\mathbb{E}\left\vert \left\langle \mathbf{u}%
_{m^{\prime}},\mathbf{\varepsilon}\right\rangle _{N}\right\vert ^{p}%
\leq\mathbb{E}\left\Vert \mathbf{\varepsilon}_{1}\right\Vert ^{p}N^{\frac
{p}{2}}+\mathbb{E}\left\Vert \mathbf{\varepsilon}_{1}\right\Vert ^{p}%
N^{\frac{p}{2}}.
\]

Then for some constant $c^{\prime}\left(  p\right)  $ that only depends on $p$%
\[
\mathbb{E}\left\vert \left\langle \mathbf{u}_{m^{\prime}},\mathbf{\varepsilon
}\right\rangle _{N}\right\vert ^{p}\leq c^{\prime}\left(  p\right)
\mathbb{E}\left\Vert \mathbf{\varepsilon}_{1}\right\Vert ^{p}N^{-\frac{p}{2}}.
\]

By this last inequality and (\ref{eq17}) we get that%
\begin{equation}
\mathbb{P}\left(  \left\vert \left\langle \mathbf{u}_{m^{\prime}%
},\mathbf{\varepsilon}\right\rangle _{N}\right\vert \geq t\right)  \leq
c^{\prime}\left(  p\right)  \mathbb{E}\left\Vert \mathbf{\varepsilon}%
_{1}\right\Vert ^{p}N^{-\frac{p}{2}}t^{-p}.\label{eq21}%
\end{equation}

Let $\upsilon$ be some positive number depending on $\eta$ only to be chosen
later. We take $t$ such that $Nt^{2}=\min\left(  \upsilon,\frac{\alpha}%
{3}\right)  \left(  L_{m^{\prime}}D_{m^{\prime}}+x\right)  \delta_{m^{\prime}%
}^{2}$ and set $N\widetilde{p}_{2}\left(  m^{\prime}\right)  =\upsilon
L_{m^{\prime}}D_{m^{\prime}}\delta_{m^{\prime}}^{2}$. We get%
\begin{align}
P_{2,m^{\prime}}\left(  x\right)   & =\mathbb{P}\left(  \frac{1}{\alpha
}\left(  \left\langle \mathbf{u}_{m^{\prime}},\mathbf{\varepsilon
}\right\rangle _{N}^{2}-\widetilde{p}_{2}\left(  m^{\prime}\right)  \right)
\geq\frac{x\delta_{m^{\prime}}^{2}}{3N}\right) \nonumber\\
& =\mathbb{P}\left(  N\left\langle \mathbf{u}_{m^{\prime}},\mathbf{\varepsilon
}\right\rangle _{N}^{2}\geq N\widetilde{p}_{2}\left(  m^{\prime}\right)
+\alpha\frac{\delta_{m^{\prime}}^{2}}{3}x\right) \nonumber\\
& =\mathbb{P}\left(  N\left\langle \mathbf{u}_{m^{\prime}},\mathbf{\varepsilon
}\right\rangle _{N}^{2}\geq\upsilon L_{m^{\prime}}D_{m^{\prime}}%
\delta_{m^{\prime}}^{2}+\alpha\frac{\delta_{m^{\prime}}^{2}}{3}x\right)
\nonumber\\
& \leq\mathbb{P}\left( \left\vert \left\langle \mathbf{u}_{m^{\prime}}%
,\mathbf{\varepsilon}\right\rangle _{N}\right\vert\geq N^{-\frac{1}{2}}\sqrt
{\min\left(  \upsilon,\frac{\alpha}{3}\right)  }\sqrt{\left(  L_{m^{\prime}%
}D_{m^{\prime}}+x\right)  }\delta_{m^{\prime}}\right) \nonumber\\
& \leq c^{\prime}\left(  p\right)  \mathbb{E}\left\Vert \mathbf{\varepsilon
}_{1}\right\Vert ^{p}N^{-\frac{p}{2}}\frac{N^{\frac{p}{2}}}{\left(
\min\left(  \upsilon,\frac{\alpha}{3}\right)  \right)  ^{\frac{p}{2}}\left(
L_{m^{\prime}}D_{m^{\prime}}+x\right)  ^{\frac{p}{2}}\delta_{m}^{p}%
}\nonumber\\
& =c^{\prime\prime}\left(  p,\eta\right)  \frac{\mathbb{E}\left\Vert
\mathbf{\varepsilon}_{1}\right\Vert ^{p}}{\delta_{m}^{p}}\frac{1}{\left(
L_{m^{\prime}}D_{m^{\prime}}+x\right)  ^{\frac{p}{2}}}.\label{eq22}%
\end{align}
The last inequality holds using (\ref{eq21}).

We now bound $P_{1,m^{\prime}}\left(  x\right)  $ for those $m^{\prime}\in
\mathcal{M}$ such that $D_{m^{\prime}}\geq1$. By using our version of Corollary 5.1
in Baraud with $\widetilde{A}=$ $\mathbf{P}_{m^{\prime}}$, ${\rm Tr}\left(
\widetilde{A}\right)  =D_{m^{\prime}}$ and $\rho\left(  \widetilde{A}\right)
=1$, we obtain from (\ref{BoundCor5.1Baraud}) that for any positive
$x_{m^{\prime}}$
\begin{equation}
\mathbb{P}\left(  N\left\Vert \mathbf{P}_{m^{\prime}}\mathbf{\varepsilon
}\right\Vert _{N}^{2}\geq\delta_{m^{\prime}}^{2}D_{m^{\prime}}+2\delta
_{m^{\prime}}^{2}\sqrt{D_{m^{\prime}}x_{m^{\prime}}}+\delta_{m^{\prime}}%
^{2}D_{m^{\prime}}x_{m^{\prime}}\right)  \leq C\left(  p\right)
\frac{\mathbb{E}\left\Vert \varepsilon_{1}\right\Vert ^{p}}{\delta_{m^{\prime
}}^{p}}D_{m^{\prime}}x_{m^{\prime}}^{-\frac{p}{2}}.\label{eq23}%
\end{equation}

Since for any $\beta>0$, $2\sqrt{D_{m^{\prime}}x_{m^{\prime}}}\leq\beta
D_{m^{\prime}}+\beta^{-1}x_{m^{\prime}}$ then (\ref{eq23}) imply that%
\begin{equation}
\mathbb{P}\left(  N\left\Vert \mathbf{P}_{m^{\prime}}\mathbf{\varepsilon
}\right\Vert _{N}^{2}\geq\left(  1+\beta\right)  D_{m^{\prime}}\delta
_{m^{\prime}}^{2}+\left(  1+\beta^{-1}\right)  x_{m^{\prime}}\delta
_{m^{\prime}}^{2}\right)  \leq C\left(  p\right)  \frac{\mathbb{E}\left\Vert
\varepsilon_{1}\right\Vert ^{p}}{\delta_{m^{\prime}}^{p}}D_{m^{\prime}%
}x_{m^{\prime}}^{-\frac{p}{2}}.\label{eq24}%
\end{equation}

Now for some number $\beta$ depending on $\eta$ only to be chosen later, we
take $x_{m^{\prime}}=\left(  1+\beta^{-1}\right)  \min\left(  \upsilon
,\frac{\left(  2-\alpha\right)  ^{-1}}{3}\right)  \left(  L_{m^{\prime}%
}D_{m^{\prime}}+x\right)  $ and $N\widetilde{p}_{1}\left(  m^{\prime}\right)
=\upsilon L_{m^{\prime}}D_{m^{\prime}}\delta_{m^{\prime}}^{2}+\left(
1+\beta\right)  D_{m^{\prime}}\delta_{m^{\prime}}^{2}$. By (\ref{eq24}) this
gives%
\begin{align}
P_{1,m^{\prime}}\left(  x\right)   & =\mathbb{P}\left(  \left\Vert
\mathbf{P}_{m^{\prime}}\mathbf{\varepsilon}\right\Vert _{N}^{2}-\widetilde
{p}_{1}\left(  m^{\prime}\right)  \geq\frac{\left(  2-\alpha\right)  ^{-1}x\delta_{m^{\prime}}^{2}}%
{3N}\right) \nonumber\\
& =\mathbb{P}\left(  N\left\Vert \mathbf{P}_{m^{\prime}}\mathbf{\varepsilon
}\right\Vert _{N}^{2}\geq\upsilon L_{m^{\prime}}D_{m^{\prime}}\delta
_{m^{\prime}}^{2}+\left(  1+\beta\right)  D_{m^{\prime}}\delta_{m^{\prime}%
}^{2}+\frac{\left(  2-\alpha\right)  ^{-1}}{3}x\delta_{m^{\prime}}^{2}\right)
\nonumber\\
& \leq\mathbb{P}\left(  N\left\Vert \mathbf{P}_{m^{\prime}}\mathbf{\varepsilon
}\right\Vert _{N}^{2}\geq\left(  1+\beta\right)  D_{m^{\prime}}\delta
_{m^{\prime}}^{2}+\left(  1+\beta^{-1}\right)  x_{m^{\prime}}\delta
_{m^{\prime}}^{2}\right) \nonumber\\
& \leq c\left(  p\right)  \frac{\mathbb{E}\left\Vert \varepsilon
_{1}\right\Vert ^{p}}{\delta_{m^{\prime}}^{p}}D_{m^{\prime}}x_{m^{\prime}%
}^{-\frac{p}{2}}\leq c^{\prime}\left(  p,\eta\right)  \frac{\mathbb{E}%
\left\Vert \varepsilon_{1}\right\Vert ^{p}}{\delta_{m^{\prime}}^{p}}%
\frac{D_{m^{\prime}}}{\left(  L_{m^{\prime}}D_{m^{\prime}}+x\right)
^{\frac{p}{2}}}.\label{eq25}%
\end{align}

Gathering (\ref{eq22}), (\ref{eq25}) and (\ref{eq16}) we get that
\begin{align*}
\mathbb{P}\left(  \mathcal{H}\left(  \mathbf{f}\right)  \geq\frac
{x\delta_{m^{\prime}}^{2}}{N\left(  1-\alpha\right)  }\right)   & \leq
\sum\limits_{m^{\prime}\in \mathcal{M}}P_{1,m^{\prime}}\left(  x\right)
+\sum\limits_{m^{\prime}\in \mathcal{M}}P_{2,m^{\prime}}\left(  x\right) \\
& \leq \sum\limits_{m^{\prime}\in \mathcal{M}}c^{\prime}\left(  p,\eta\right)
\frac{\mathbb{E}\left\Vert \varepsilon_{1}\right\Vert ^{p}}{\delta_{m^{\prime
}}^{p}}\frac{D_{m^{\prime}}}{\left(  L_{m^{\prime}}D_{m^{\prime}}+x\right)
^{\frac{p}{2}}} \\
+ & \sum\limits_{m^{\prime}\in \mathcal{M}}c^{\prime\prime}\left(
p,\eta\right)  \frac{\mathbb{E}\left\Vert \varepsilon_{1}\right\Vert ^{p}%
}{\delta_{m^{\prime}}^{p}}\frac{1}{\left(  L_{m^{\prime}}D_{m^{\prime}%
}+x\right)  ^{\frac{p}{2}}}.
\end{align*}

Since $\frac{1}{\left(  1-\alpha\right)  }=\left(  1+2\eta^{-1}\right)  $,
then (\ref{ineqPH}) holds:
\begin{align*}
\mathbb{P}\left(  \mathcal{H}\left(  \mathbf{f}\right)  \geq\left(
1+2\eta^{-1}\right)  \frac{x\delta_{m^{\prime}}^{2}}{N}\right)   &  \leq
\sum\limits_{m^{\prime}\in \mathcal{M}}\frac{\mathbb{E}\left\Vert \varepsilon
_{1}\right\Vert ^{p}}{\delta_{m^{\prime}}^{p}\left(  L_{m^{\prime}%
}D_{m^{\prime}}+x\right)  ^{\frac{p}{2}}}\max\left(  D_{m^{\prime}},1\right)
\left(  c^{\prime}\left(  p,\eta\right)  +c^{\prime\prime}\left(
p,\eta\right)  \right)  \\
&  =c\left(  p,\eta\right)  \frac{\mathbb{E}\left\Vert \varepsilon
_{1}\right\Vert ^{p}}{\delta_{m^{\prime}}^{p}}\sum\limits_{m^{\prime}\in
\mathcal{M}}\frac{D_{m^{\prime}}\vee1}{\left(  L_{m^{\prime}}D_{m^{\prime}%
}+x\right)  ^{\frac{p}{2}}}.
\end{align*}

It remains to choose $\beta$ and $\delta$ for (\ref{eq14}) to hold (we recall
that $\alpha=\frac{2}{2+\eta}$). This is the case if $\left(  2-\alpha\right)
\left(  1+\beta\right)  =1+\eta$ and $\left(  2-\alpha+\alpha^{-1}\right)
\delta=1$, therefore we take $\beta=\frac{\eta}{2}$ and $\delta=\left[
1+\frac{\eta}{2}+2\frac{(1+\eta)}{(2+\eta)}\right]  ^{-1}$.
\end{proof}

\subsection{Proof of the concentration inequality}

Proof of Proposition~\eqref{ExtCor5.1Baraud}
\begin{proof}
Denote by $\tau^{2}$ the following expression:%
\[
\tau^{2}:=\mathbb{E}\left\Vert \mathbf{P}_{m}\mathbf{\varepsilon}\right\Vert
^{2}=\mathbb{E}\left(  \mathbf{\varepsilon}^{T}\mathbf{P}_{m}%
\mathbf{\varepsilon}\right)  ={\rm Tr}\left(  \mathbf{P}_{m}\left(  I_{N}\otimes
\Phi\right)  \right)  .
\]

Then we have that%
\begin{align*}
\tau^{2}  & ={\rm Tr}\left(  \mathbf{P}_{m}\left(  I_{N}\otimes\Phi\right)
\mathbf{P}_{m}\right)  \leq\lambda_{\max}\left(  I_{N}\otimes\Phi\right)
{\rm Tr}\left(  \mathbf{P}_{m}^{2}\right)  =\lambda_{\max}\left(  I_{N}\otimes
\Phi\right)  {\rm Tr}\left(  \mathbf{P}_{m}\right) \\
& =\lambda_{\max}\left(  \Phi\right)  {\rm Tr}\left(  \mathbf{P}_{m}\right)
=\lambda_{\max}\left(  \Phi\right)  D_{m}.
\end{align*}

We have that $\eta^{2}\left(  \mathbf{\varepsilon}\right)
:=\mathbf{\varepsilon}^{T}\widetilde{A}\mathbf{\varepsilon}$, where
$\widetilde{A}=A^{T}A$. Then%
\begin{align*}
\eta^{2}\left(  \mathbf{\varepsilon}\right)    & =\left\Vert
A\mathbf{\varepsilon}\right\Vert ^{2}=\left[  \underset{\left\Vert
u\right\Vert \leq1}{\sup}\left\langle A\mathbf{\varepsilon},\mathbf{u}%
\right\rangle \right]  ^{2}=\left[  \underset{\left\Vert u\right\Vert \leq
1}{\sup}\sum\limits_{i=1}^{Nd}\left(  A\mathbf{\varepsilon}\right)
_{i}\mathbf{u}_{i}\right]  ^{2}\\
& =\left[  \underset{\left\Vert \mathbf{u}\right\Vert \leq1}{\sup}\left\langle
\mathbf{\varepsilon},A^{T}\mathbf{u}\right\rangle \right]  ^{2}=\left[
\underset{\left\Vert \mathbf{u}\right\Vert \leq1}{\sup}\sum\limits_{i=1}%
^{N}\left\langle \mathbf{\varepsilon}_{i},\left(  A^{T}\mathbf{u}\right)
_{i}\right\rangle \right]  ^{2}\\
& =\left[  \underset{\left\Vert \mathbf{u}\right\Vert \leq1}{\sup}%
\sum\limits_{i=1}^{N}\left\langle \mathbf{\varepsilon}_{i},A_{i}^{T}%
\mathbf{u}\right\rangle \right]  ^{2}=\left[  \underset{\left\Vert
\mathbf{u}\right\Vert \leq1}{\sup}\sum\limits_{i=1}^{N}\sum\limits_{j=1}%
^{d}\varepsilon_{ij}\left(  A_{i}^{T}\mathbf{u}\right)  _{j}\right]  ^{2}%
\end{align*}
with $A=\left(  A_{1}\mid...\mid A_{N}\right)  $, where $A_{i}$ is a $\left(
Nd\right)  \times d$ matrix.

Now take $\mathcal{G}= \{  g_{\mathbf{u}}:g_{\mathbf{u}}\left(
\mathbf{x}\right)  =\sum\limits_{i=1}^{N}\left\langle \mathbf{x}_{i},A_{i}%
^{T}\mathbf{u}\right\rangle =\sum\limits_{i=1}^{N}\left\langle B_{i}%
\mathbf{x},B_{i}A^{T}\mathbf{u}\right\rangle ,\;\mathbf{u},\mathbf{x}=\left(
\mathbf{x}_{1},\dots,\mathbf{x}_{N} \right)^{'}  \in\mathbb{R}^{\left(  Nd\right)  },\;\left\Vert \mathbf{u}%
\right\Vert \leq1 \}  .$

Let $M_{i}=\left[\mathbf{0},\dots,
\mathbf{0},I_{d},\mathbf{0},\dots,
\mathbf{0}
\right]^{'}  \in\mathbb{R}^{\left(  Nd\right)  \times\left(  Nd\right)  }$, where
$I_{d}$ is the $i-$th block of $M_{i}$, $B_{i}=\left[  0,...,0,I_{d}%
,0,...0\right]  \in\mathbb{R}^{\left(  Nd\right)  \times\left(  Nd\right)  }$,
$\mathbf{\varepsilon}_{i}=B_{i}$ $\mathbf{\varepsilon}$ and $M_{i}$
$\mathbf{\varepsilon}=\left[
\mathbf{0},\dots,
\mathbf{0},
\mathbf{\varepsilon}_{i},
\mathbf{0},\dots,
\mathbf{0}
\right]^{'}  $.

Then%
\[
\eta\left(  \mathbf{\varepsilon}\right)  =\underset{\left\Vert \mathbf{u}%
\right\Vert \leq1}{\sup}\sum\limits_{i=1}^{N}g_{\mathbf{u}}\left(
M_{i}\mathbf{\varepsilon}\right)  .
\]

Now take $\mathbf{U}_{i}=M_{i}$ $\mathbf{\varepsilon}$, $\mathbf{\varepsilon
}\in\mathbb{R}^{\left(  Nd\right)  }$. Then for each positive number $t$ and
$p>0$%
\begin{align}
\mathbb{P}\left(  \eta\left(  \mathbf{\varepsilon}\right)  \geq\mathbb{E}%
\left(  \eta\left(  \mathbf{\varepsilon}\right)  \right)  +t\right)   &
\leq\mathbb{P}\left(  \left\vert \eta\left(  \mathbf{\varepsilon}\right)
-\mathbb{E}\left(  \eta\left(  \mathbf{\varepsilon}\right)  \right)
\right\vert >t\right) \nonumber\\
& \leq t^{-p}\mathbb{E}\left(  \left\vert \eta\left(  \mathbf{\varepsilon
}\right)  -\mathbb{E}\left(  \eta\left(  \mathbf{\varepsilon}\right)  \right)
\right\vert ^{p}\right)  \text{ by Markov inequality}\nonumber\\
 \leq c\left(  p\right)  t^{-p} & \left\{  \mathbb{E}\left(  \underset
{i=1,...,N}{\max}\underset{\left\Vert \mathbf{u}\right\Vert \leq1}{\sup
}\left\vert \left\langle \mathbf{\varepsilon}_{i},A_{i}^{T}\mathbf{u}%
\right\rangle \right\vert ^{p}\right)  +\left[  \mathbb{E}\left(
\underset{\left\Vert \mathbf{u}\right\Vert \leq1}{\sup}\sum\limits_{i=1}%
^{N}\left(  \left\langle \mathbf{\varepsilon}_{i},A_{i}^{T}\mathbf{u}%
\right\rangle \right)  ^{2}\right)  \right]  ^{p/2}\right\} \nonumber\\
& =c\left(  p\right)  t^{-p}\left(  \mathbb{E}_{1}+\mathbb{E}_{2}%
^{p/2}\right)  .\label{E1+E2}%
\end{align}

We start by bounding $\mathbb{E}_{1}$. For all $\mathbf{u}$ such that
$\left\Vert \mathbf{u}\right\Vert \leq1$ and $i\in\left\{  1,...,N\right\}
$,
\[
\left\Vert A_{i}^{T}\mathbf{u}\right\Vert ^{2}\leq\left\Vert A^{T}%
\mathbf{u}\right\Vert ^{2}\leq\rho^{2}\left(  A\right)  ,
\]
where $\rho\left(  M\right)  =\underset{x\neq0}{\sup}\frac{\left\Vert
Mx\right\Vert }{\left\Vert x\right\Vert }\;$for all matrix $M$. For $p\geq2$
we have that $\left\Vert A_{i}\mathbf{u}\right\Vert ^{p}\leq\rho^{p-2}\left(
A\right)  \left\Vert A_{i}\mathbf{u}\right\Vert ^{2}$, then%
\[
\left\vert \left\langle \mathbf{\varepsilon}_{i},A_{i}^{T}\mathbf{u}%
\right\rangle \right\vert ^{p}\leq\left[  \left\Vert \mathbf{\varepsilon}%
_{i}\right\Vert \left\Vert A_{i}^{T}\mathbf{u}\right\Vert \right]  ^{p}%
\leq\rho^{p-2}\left(  A\right)  \left\Vert \mathbf{\varepsilon}_{i}\right\Vert
^{p}\left\Vert A_{i}^{T}\mathbf{u}\right\Vert ^{p}.
\]

Therefore%
\[
\mathbb{E}_{1}\leq\rho^{p-2}\left(  A\right)  \mathbb{E}\left(  \underset
{\left\Vert \mathbf{u}\right\Vert =1}{\sup}\sum\limits_{i=1}^{N}\left\Vert
\mathbf{\varepsilon}_{i}\right\Vert ^{p}\left\Vert A_{i}^{T}\mathbf{u}%
\right\Vert ^{2}\right)  .
\]

Since $\left\Vert \mathbf{u}\right\Vert \leq1$, $\forall i=1,...,N$
\begin{align*}
\left\Vert A_{i}^{T}\mathbf{u}\right\Vert ^{2}  & =\mathbf{u}^{T}A_{i}%
A_{i}^{T}\mathbf{u}\leq\rho\left(  A_{i}A_{i}^{T}\right) \\
& \leq {\rm Tr}\left(  A_{i}A_{i}^{T}\right)  ,
\end{align*}
then%
\[
\sum\limits_{i=1}^{N}\left\Vert A_{i}^{T}\mathbf{u}\right\Vert ^{2}\leq
\sum\limits_{i=1}^{N}{\rm Tr}\left(  A_{i}A_{i}^{T}\right)  ={\rm Tr}\left(
\sum\limits_{i=1}^{N}A_{i}A_{i}^{T}\right)  ={\rm Tr}\left(  \widetilde{A}\right)  .
\]
Thus,%
\begin{equation}
\mathbb{E}_{1}\leq\rho^{p-2}\left(  A\right)  {\rm Tr}\left(  \widetilde{A}\right)
\mathbb{E}\left(  \left\Vert \mathbf{\varepsilon}_{i}\right\Vert ^{p}\right)
.\label{ineqE1}%
\end{equation}

We now bound $\mathbb{E}_{2}$ via a truncation argument. Since for all
$\mathbf{u}$ such that $\left\Vert \mathbf{u}\right\Vert \leq1$ and
$i\in\left\{  1,...,N\right\}  $, $\left\Vert A^{T}\mathbf{u}\right\Vert
^{2}\leq\rho^{2}\left(  A\right)  $, for any positive number $c$ to be
specified later we have that%
\begin{align}
\mathbb{E}_{2}  & \leq\mathbb{E}\left(  \underset{\left\Vert \mathbf{u}%
\right\Vert \leq1}{\sup}\sum\limits_{i=1}^{N}\left\Vert \mathbf{\varepsilon
}_{i}\right\Vert ^{2}\left\Vert A_{i}^{T}\mathbf{u}\right\Vert ^{2}1_{\left\{
\left\Vert \mathbf{\varepsilon}_{i}\right\Vert \leq c\right\}  }\right)
+\mathbb{E}\left(  \underset{\left\Vert \mathbf{u}\right\Vert \leq1}{\sup}%
\sum\limits_{i=1}^{N}\left\Vert \mathbf{\varepsilon}_{i}\right\Vert
^{2}\left\Vert A_{i}^{T}\mathbf{u}\right\Vert ^{2}1_{\left\{  \left\Vert
\mathbf{\varepsilon}_{i}\right\Vert >c\right\}  }\right) \nonumber\\
& \leq\mathbb{E}\left(  c^{2}\underset{\left\Vert \mathbf{u}\right\Vert \leq
1}{\sup}\sum\limits_{i=1}^{N}\left\Vert A_{i}^{T}\mathbf{u}\right\Vert
^{2}1_{\left\{  \left\Vert \mathbf{\varepsilon}_{i}\right\Vert \leq c\right\}
}\right)  +\mathbb{E}\left(  \underset{\left\Vert \mathbf{u}\right\Vert \leq
1}{\sup}\sum\limits_{i=1}^{N}\left\Vert \mathbf{\varepsilon}_{i}\right\Vert
^{2}\left\Vert A_{i}^{T}\mathbf{u}\right\Vert ^{2}1_{\left\{  \left\Vert
\mathbf{\varepsilon}_{i}\right\Vert >c\right\}  }\right) \nonumber\\
& \leq c^{2}\rho^{2}\left(  A\right)  +c^{2-p}\mathbb{E}\left(  \underset
{\left\Vert \mathbf{u}\right\Vert \leq1}{\sup}\sum\limits_{i=1}^{N}\left\Vert
A_{i}\mathbf{u}\right\Vert ^{2}\left\Vert \mathbf{\varepsilon}_{i}\right\Vert
^{p}\right)  \text{ }\nonumber\\
& \leq c^{2}\rho^{2}\left(  A\right)  +c^{2-p}\mathbb{E}\left(  \left\Vert
\mathbf{\varepsilon}_{i}\right\Vert ^{p}\right)  {\rm Tr}\left(  \widetilde
{A}\right) \label{ineqE2}%
\end{align}
using the bound obtained for $\mathbb{E}_{1}.$ It remains to take
$c^{p}=\mathbb{E}\left(  \left\Vert \mathbf{\varepsilon}_{i}\right\Vert
^{p}\right)  {\rm Tr}\left(  \widetilde{A}\right)  /\rho^{2}\left(  A\right)  $ to
get that:%
\[
\mathbb{E}_{2}\leq c^{2}\rho^{2}\left(  A\right)  +c^{2}\rho^{2}\left(
A\right)  =2c^{2}\rho^{2}\left(  A\right)  ,
\]
therefore%
\begin{equation}
\mathbb{E}_{2}^{p/2}\leq 2^{p/2}c^{p}\rho^{p}\left(  A\right)
,\label{ineqE2^p/2}%
\end{equation}
which implies that%
\[
2^{-p/2}\mathbb{E}_{2}^{p/2}\leq\mathbb{E}\left(  \left\Vert
\mathbf{\varepsilon}_{1}\right\Vert ^{p}\right)  {\rm Tr}\left(  \widetilde
{A}\right)  \rho^{p-2}\left(  A\right)  .
\]

We straightforwardly derive from (\ref{E1+E2}) that%
\[
\mathbb{P}\left(  \eta^{2}\left(  \mathbf{\varepsilon}\right)  \geq\left[
\mathbb{E}\left(  \eta\left(  \mathbf{\varepsilon}\right)  \right)  \right]
^{2}+2\mathbb{E}\left(  \eta\left(  \mathbf{\varepsilon}\right)  \right)
t+t^{2}\right)  \leq c\left(  p\right)  t^{-p}\left(  \mathbb{E}%
_{1}+\mathbb{E}_{2}^{p/2}\right)  .
\]

Since $\left[  \mathbb{E}\left(  \eta\left(  \mathbf{\varepsilon}\right)
\right)  \right]  ^{2}\leq\mathbb{E}\left(  \eta^{2}\left(
\mathbf{\varepsilon}\right)  \right)  $, (\ref{ineqE1}) and (\ref{ineqE2^p/2})
imply that%
\begin{align}
& \mathbb{P}\left(  \eta^{2}\left(  \mathbf{\varepsilon}\right)  \geq
\mathbb{E}\left(  \eta^{2}\left(  \mathbf{\varepsilon}\right)  \right)
+2\sqrt{\mathbb{E}\left(  \eta^{2}\left(  \mathbf{\varepsilon}\right)
\right)  t^{2}}+t^{2}\right)    \leq c\left(  p\right)  t^{-p}\left(
\mathbb{E}_{1}+\mathbb{E}_{2}^{p/2}\right) \nonumber\\
 & \leq c\left(  p\right)  t^{-p}   \left(  \rho^{p-2} \left(  A\right)  {\rm Tr}\left(
\widetilde{A}\right)  \mathbb{E}\left(  \left\Vert \mathbf{\varepsilon}%
_{i}\right\Vert ^{p}\right)  +2^{p/2}\mathbb{E}\left(  \left\Vert
\mathbf{\varepsilon}_{1}\right\Vert ^{p}\right)  {\rm Tr} \left(  \widetilde
{A}\right)  \rho^{p-2}\left(  A\right)  \right) \nonumber\\
& \leq c^{\prime}\left(  p\right)  t^{-p}\rho^{p-2}\left(  A\right)  {\rm Tr}\left(
\widetilde{A}\right)  \mathbb{E}\left(  \left\Vert \mathbf{\varepsilon}%
_{i}\right\Vert ^{p}\right)  ,\label{ineq5}%
\end{align}
for all $t>0$. Moreover%
\begin{align*}
\mathbb{E}\left(  \eta^{2}\left(  \mathbf{\varepsilon}\right)  \right)   &
=\mathbb{E}\left(  \mathbf{\varepsilon}^{T}\widetilde{A}\mathbf{\varepsilon
}\right)  =\mathbb{E}\left(  \left\Vert A\mathbf{\varepsilon}\right\Vert
^{2}\right)  =\mathbb{E}\left(  \sum\limits_{i=1}^{N}\left\Vert A_{i}%
\mathbf{\varepsilon}_{i}\right\Vert ^{2}\right) \\
& =\sum\limits_{i=1}^{N}\mathbb{E}\left(  {\rm Tr}\mathbf{\varepsilon}_{i}^{T}%
A_{i}^{T}A_{i}\mathbf{\varepsilon}_{i}\right)  =\sum\limits_{i=1}^{N}%
{\rm Tr}A_{i}^{T}A_{i}\mathbb{E}\left(  \mathbf{\varepsilon}_{i}\mathbf{\varepsilon
}_{i}^{T}\right) \\
& ={\rm Tr}\left(  \sum\limits_{i=1}^{N}A_{i}^{T}A_{i}\right)  \Phi.
\end{align*}

But it is better to use that%
\begin{align}
\mathbb{E}\left(  \eta^{2}\left(  \mathbf{\varepsilon}\right)  \right)   &
={\rm Tr} \left( \widetilde{A}\mathbf{\varepsilon\varepsilon}^{T}\right)={\rm Tr} \left( \widetilde{A}\left(
I_{N}\otimes\Phi\right)\right)  ={\rm Tr}\left(A^{T}A\left(  I_{N}\otimes\Phi\right)\right)  ={\rm Tr}\left( A\left(
I_{N}\otimes\Phi\right)  A^{T}\right) \nonumber\\
& \leq\lambda_{\max}\left(  I_{N}\otimes\Phi\right)  {\rm Tr}\left(AA^{T}\right)=\lambda_{\max
}\left(  I_{N}\otimes\Phi\right)  {\rm Tr}\left(\widetilde{A}\right)=\lambda_{\max}\left(
Q\right)  {\rm Tr}\left(\widetilde{A}\right),\label{ineq6}%
\end{align}
for $Q=I_{N}\otimes\Phi$.

Using (\ref{ineq6}), take $t^{2}=\rho\left(  \widetilde{A}\right)
\lambda_{\max}\left(  I_{N}\otimes\Phi\right)  x>0$ in (\ref{ineq5}) to get
that
\[
\mathbb{P}\left(  \eta^{2}\left(  \mathbf{\varepsilon}\right)  \geq
\lambda_{\max}\left(  Q\right)  {\rm Tr}\left(  \widetilde{A}\right)  +2\sqrt
{\lambda_{\max}\left(  Q\right)  {\rm Tr}\left(  \widetilde{A}\right)  \rho\left(
\widetilde{A}\right)  \lambda_{\max}\left(  Q\right)  x}+\rho\left(
\widetilde{A}\right)  \lambda_{\max}\left(  Q\right)  x\right)  \]
\[ \leq
c^{\prime}\left(  p\right)  \rho^{-p/2}\left(  \widetilde{A}\right)
\lambda_{\max}^{-p/2}\left(  Q\right)  x^{-p/2}\rho^{p-2}\left(  A\right)
{\rm Tr}\left(  \widetilde{A}\right)  \mathbb{E}\left(  \left\Vert
\mathbf{\varepsilon}_{i}\right\Vert ^{p}\right)  .
\]

Since $\rho\left(  \widetilde{A}\right)  =\rho^{2}\left(  A\right)  $ (with
the Euclidean norm) the desired result follows:%
\begin{equation*}
\mathbb{P}\left(  \eta^{2}\left(  \mathbf{\varepsilon}\right)  \geq
\lambda_{\max}\left(  Q\right)  {\rm Tr}\left(  \widetilde{A}\right)  +2\lambda
_{\max}\left(  Q\right)  \sqrt{\rho\left(  \widetilde{A}\right)  {\rm Tr}\left(
\widetilde{A}\right)  x}+\lambda_{\max}\left(  Q\right)  \rho\left(
\widetilde{A}\right)  x\right)
\end{equation*}
\begin{equation}
\leq c^{\prime}\left(  p\right)
\frac{\mathbb{E}\left\Vert \mathbf{\varepsilon}_{1}\right\Vert ^{p}}{\left(
\sqrt{\lambda_{\max}\left(  Q\right)  }\right)  ^{p}}\frac{{\rm Tr}\left(
\widetilde{A}\right)  }{\rho\left(  \widetilde{A}\right)  x^{p/2}%
}.\label{ineq7}%
\end{equation}
\end{proof}
\bibliographystyle{alpha}
\bibliography{CovSelect2}
\noindent J. {\sc Bigot} \& J-M. {\sc Loubes}\hfill  R. Biscay \& L. Mu\~niz\\
\noindent Equipe de probabilit\'es et statistique,\hfill Instituto de Cibernetica, Matematica y Fisica,\\ 
Institut de Math\'{e}matique de Toulouse,\hfill  Departamiento de Matematicas,\\
UMR5219, Universit\'e de Toulouse,\hfill Universidad Central de la Havana,\\
\noindent 31000 Toulouse  {\sc France}\hfill Ciudad Havana  {\sc Cuba}\\
\end{document}

%% file: macros.tex
\newcommand{\DOIFPDF}[2]{\ifx\pdfoutput\undefined #2\else#1\fi}

\newcommand{\EM}{\ensuremath}

\newcommand{\PDFABLE}[2]{%
\newif\ifpdf%
\ifx\pdfoutput\undefined\pdffalse%
\else\pdftrue\pdfoutput=1\pdfcompresslevel=9\fi%
\ifpdf%
 \usepackage#1
 \usepackage[pdftex,%
             a4paper,%
             colorlinks,%
             citecolor=blue,%
             pagebackref,%
             plainpages=false]{hyperref}%
\else%
 \usepackage#2
 \usepackage{url}
\fi}

\makeatletter
\newcommand{\@THMSTYLES}{%
  \newtheoremstyle{bodyrm}
  {3pt}
  {3pt}
  {}
  {}
  {\bfseries\sffamily}
  {.}
  { }
  {}
  \newtheoremstyle{bodyit}
  {3pt}
  {3pt}
  {\itshape}
  {}
  {\bfseries\sffamily}
  {.}
  { }
  {}
}
\newcommand{\THMEN}{%
  \@THMSTYLES
  \theoremstyle{bodyit}
  \newtheorem{thm}{Theorem}[section]%
  \newtheorem{cor}[thm]{Corollary}%
  \newtheorem{prop}[thm]{Proposition}%
  \newtheorem{lem}[thm]{Lemma}%
  \theoremstyle{bodyrm}%
  \newtheorem{defi}[thm]{Definition}%
  \newtheorem{xpl}[thm]{Example}%
  \newtheorem{exo}[thm]{Exercise}%
  \newtheorem{hyp}[thm]{Hypothesis}%
  \newtheorem{eur}[thm]{Heuristics}%
  \newtheorem{pro}[thm]{Problem}%
  \newtheorem{rem}[thm]{Remark}%
  \newtheorem{prp}[thm]{Property}%
}
\newcommand{\THMFR}{%
  \@THMSTYLES
  \theoremstyle{bodyit}
  \newtheorem{thm}{Théorème}[section]%
  \newtheorem{cor}[thm]{Corollaire}%
  \newtheorem{prop}[thm]{Proposition}%
  \theoremstyle{bodyrm}%
  \newtheorem{rem}[thm]{Remarque}%
}
\makeatother

\makeatletter
\newcommand{\SMALLSECS}{%
 \renewcommand{\section}{\@startsection%
  {section}
  {1}
  {0em}
  {\baselineskip}
  {0.5\baselineskip}
  {\normalfont\large\bfseries}}
 \renewcommand{\subsection}{\@startsection%
  {subsection}
  {2}
  {0em}
  {\baselineskip}
  {0.25\baselineskip}
  {\normalfont\bfseries}}
}
\makeatother

\makeatletter
\providecommand{\timenow}{\@tempcnta\time
\@tempcntb\@tempcnta
\divide\@tempcntb60
\ifnum10>\@tempcntb0\fi\number\@tempcntb
\multiply\@tempcntb60
\advance\@tempcnta-\@tempcntb
:\ifnum10>\@tempcnta0\fi\number\@tempcnta}
\makeatother

\makeatletter
\newcommand{\versiondetravail}{%
 \renewcommand{\@evenfoot}{%
 \hfil{\tiny\texttt{%
   Version préliminaire, compilée le \today{} à \timenow.}\hfill}}%
 \renewcommand{\@oddfoot}{\@evenfoot}%
}
\makeatother

\newcommand{\bR}{\EM{\mathbf{R}}}

\newcommand{\Det}[1]{\mathrm{Det}\,}

\renewcommand{\leq}{\leqslant}
\renewcommand{\geq}{\geqslant}

